\documentclass[12pt,twoside, epsf]{article}

\usepackage{color,graphicx,times}

\makeatother

\let \ttorg \tt \def \tt{\ttorg \obeyspaces}

\begin{document}

\newcommand{\Across}{\raisebox{-0.25\height}{\includegraphics[width=0.5cm]{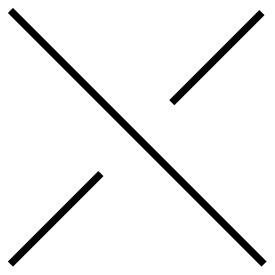}}}
\newcommand{\Asmooth}{\raisebox{-0.25\height}{\includegraphics[width=0.5cm]{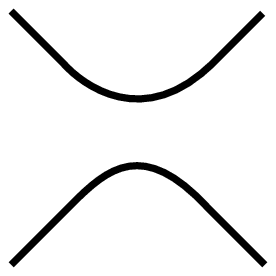}}}
\newcommand{\Bsmooth}{\raisebox{-0.25\height}{\includegraphics[width=0.5cm]{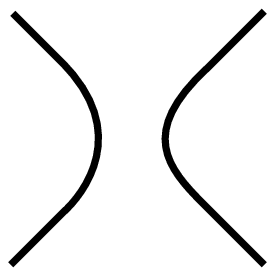}}}
\newcommand{\Rcurl}{\raisebox{-0.25\height}{\includegraphics[width=0.5cm]{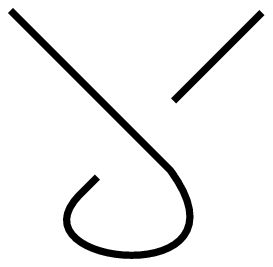}}}
\newcommand{\Lcurl}{\raisebox{-0.25\height}{\includegraphics[width=0.5cm]{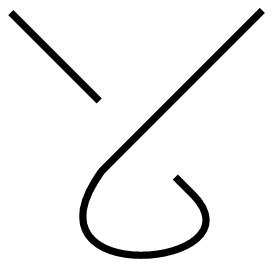}}}
\newcommand{\Arc}{\raisebox{-0.25\height}{\includegraphics[width=0.5cm]{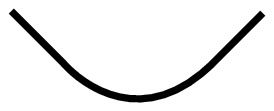}}}

\date{}

\title{\Large\bf Introduction to Virtual Knot Theory}
\author{Louis H. Kauffman\\ Department of Mathematics, Statistics \\ and Computer Science (m/c
249)    \\ 851 South Morgan Street   \\ University of Illinois at Chicago\\
Chicago, Illinois 60607-7045\\ $<$kauffman@uic.edu$>$\\}

\maketitle

\thispagestyle{empty}

\subsection*{\centering Abstract}

{\em This paper  is an introduction to virtual knot theory and an exposition of new ideas and constructions, including the parity bracket polynomial, the arrow polynomial, the parity arrow polynomial and categorifications of the arrow polynomial.}

\section{Introduction}
Virtual knot theory is an extension of classical knot theory to stabilized embeddings of circles into thickened orientable surfaces of
genus possibly greater than zero. Classical knot theory is the case of genus zero. There is a diagrammatic theory for studying virtual knots and links, 
and this diagrammatic theory lends itself to the construction of numerous new invariants of virtual knots as well as extensions of known
invariants. 
\bigbreak

This paper gives an introduction to virtual knot theory and it defines a new invariant of virtual links and flat virtual links that we call the {\it arrow polynomial} \cite{DyeKauff,ExtBr}.
 \bigbreak

One of the applications of this new invariant is to the category of flat virtual knots and links. In Section 5 we review the definitions of flat
virtuals and we prove that {\it isotopy classes of long flat virtual knots embed in the isotopy classes of 
all long virtual knots} via the ascending map $A: LVF \longrightarrow LVK$. See Section 5 for this result. This embedding of the long flat virtual knots means
that there are many invariants of them  obtained by appling any invariant $Inv$ of virtual knots via the composition with the ascending map $A.$ This
situation is in direct contrast to closed long virtual flats, where it is not so easy to define invariants. The arrow polynonmial is an
invariant of both long virtual flats and closed virtual flats.
\bigbreak  

This paper is relatively self-contained, with Section 2 reviewing the definitions of virtual knots and links, Section 3 reviewing the surface interpretations 
of virtuals and Section 4 reviewing the definitions and properties of flat virtuals. 
\bigbreak

Section 5 proves that long flat virtual knots (virtual strings in Turaev's terminology) embed into long virtual knots. This result is a very strong tool for discriminating long virtual knots from one another. We give examples and computations using this technique by applying invariants of virtual knots to corresponding ascending diagram images of flat virtual knots. A long virtual can be closed just as a long knot can be closed. It is  a remarkable property of long virtuals (both flat and not flat) that there are non-trivial long examples whose closures are trivial. Thus one would 
like to understand the kernel of the closure mapping from long virtuals to virtuals both in the flat and non-flat categories. It is hoped that the invariants discussed in this paper will further this question.
We also introduce in Section 5 the notion of the parity (odd or even) of a crossing in a virtual knot and recall the {\it odd writhe} of a virtual
knot -- an invariant that is a kind of self-linking number for virtual knots.
 \bigbreak

Section 6 is a review of the bracket polynomial and the Jones polynomial
for virtual knots and links. We recall the virtualization construction that produces infinitely many non-trivial virtual knots with unit Jones polynomial.
These examples are of interest since there are no known examples of classical non-trivial knots with unit Jones polynomial. One may conjecture that
all non-trivial examples produced by virtualization are non-classical. It may be that the virtualization of some non-trivial classical knot is 
isotopic to a classical knot, but we have no evidence that this can happen. The arrow polynomial invariant developed in this paper may help in
deciding this question of non-classicality. 
\bigbreak

Section 7 introduces the Parity Bracket, a modification of the bracket polynomial due to Vassily Manturov that takes into account the parity of crossings in a knot diagram. The Parity Bracket takes values in sums of certain reduced graphs with polynomial coefficients. We show the power of this invariant by given
an example of a virtual knot with unit Jones polynomial that is not obtained from a classical knot
by the virtualization construction of the previous section. This is the first time such an example has been
pinpointed and we do not know any method to detect such examples at this time other than
the Parity Bracket of Vassily Manturov.
\bigbreak

In Section 8 we give the definition of the arrow polynomial  ${\cal A}[K]$ and a number of examples of its calculation. These examples 
include a verification of the non-classicality of the simplest example of virtualization,
a verification that the Kishino diagram \cite{KIS} and the flat Kishino diagram are non-trivial.  The arrow polynomial is an invariant of flat diagrams by taking the specialization of its parameters so that $A=1$ and $\delta = -2.$
We give an example showing that the arrow polynomial can detect a long virtual knot whose closure is trivial. This is a capability that is beyond the reach of the Jones polynomial. As the reader will see from the definitions in section 8, the arrow polynomial is a direct generalization of the bracket model of Jones polynomial, using the oriented structure of the link diagrams. Section 8 ends with a short discussion of
how to create a parity version of the arrow polynomial. This study will be taken up in a separate 
paper.
\bigbreak 

The arrow polynomial has infinitely many variables and integer coefficients.  It has been discovered independently by Miyazawa, via a different definition in \cite{Miyazawa2,KamadaMiya}. Miyazawa does not consider the applications of this polynomial to flat virtual knots. In
\cite{ExtBr} we use the interaction of the arrow polynomial with the extended bracket polynomial to obtain results on the genus of virtual knots. In this paper  we go beyond our previous work on the arrow polynomial by extending it to long virtual knots and to examples involving categorification and parity. 
\bigbreak

In Section 9 we give a quick review of Khovanov homology, and in Section 10 we describe the work 
of Heather Dye, the Author and Vassily Manturov in categorifying the arrow polyomial. We include an
example (one of many found by Aaron Kaestner) of a pair of virtual knots that are not distinguished 
by the Jones polynomial, the bracket polynomial, the arrow polynomial  or by Khovanov homology, while the pair is distinguished by our categorification of the arrow polynomial. We also show that this
pair is distinguished by the Parity Bracket. 
\bigbreak

Section 10 ends the paper with a discussion about futher directions in virtual knot theory.
\bigbreak

\section{Virtual Knot Theory}
Knot theory
studies the embeddings of curves in three-dimensional space.  Virtual knot theory studies the  embeddings of curves in thickened surfaces of arbitrary
genus, up to the addition and removal of empty handles from the surface. Virtual knots have a special diagrammatic theory, described below,
that makes handling them
very similar to the handling of classical knot diagrams. Many structures in classical knot
theory generalize to the virtual domain.
\bigbreak  

In the diagrammatic theory of virtual knots one adds 
a {\em virtual crossing} (see Figure 1) that is neither an over-crossing
nor an under-crossing.  A virtual crossing is represented by two crossing segments with a small circle
placed around the crossing point. 
\bigbreak

Moves on virtual diagrams generalize the Reidemeister moves for classical knot and link
diagrams.  See Figure 1.  One can summarize the moves on virtual diagrams by saying that the classical crossings interact with
one another according to the usual Reidemeister moves while virtual crossings are artifacts of the attempt to draw the virtual structure in the plane. 
A segment of diagram consisting of a sequence of consecutive virtual crossings can be excised and a new connection made between the resulting
free ends. If the new connecting segment intersects the remaining diagram (transversally) then each new intersection is taken to be virtual.
Such an excision and reconnection is called a {\it detour move}.
Adding the global detour move to the Reidemeister moves completes the description of moves on virtual diagrams. In Figure 1 we illustrate a set of local
moves involving virtual crossings. The global detour move is
a consequence of  moves (B) and (C) in Figure 1. The detour move is illustrated in Figure 2.  Virtual knot and link diagrams that can be connected by a finite 
sequence of these moves are said to be {\it equivalent} or {\it virtually isotopic}.
\bigbreak

\begin{figure}[htb]
     \begin{center}
     \begin{tabular}{c}
     \includegraphics[width=10cm]{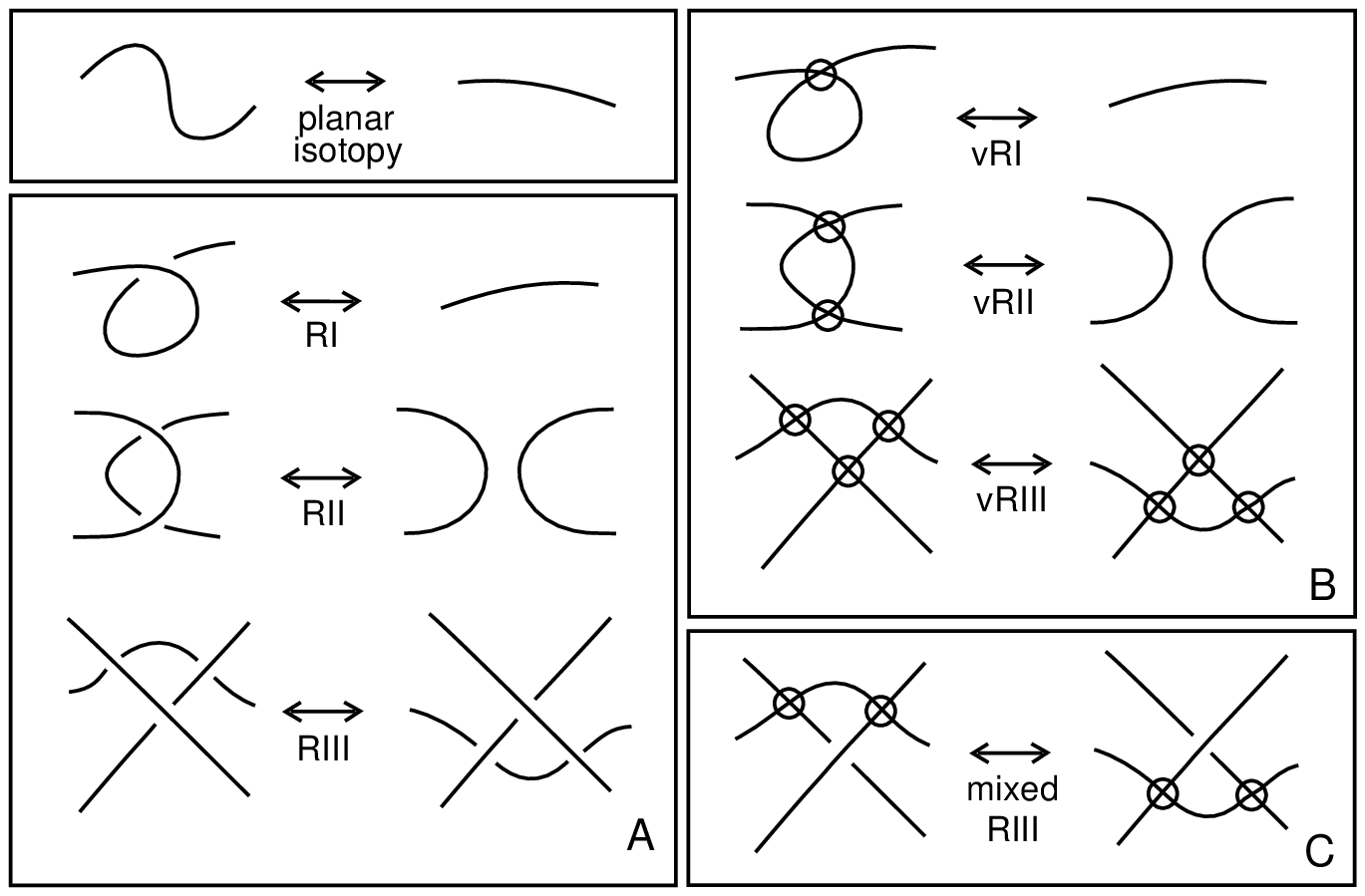}
     \end{tabular}
     \caption{\bf Moves}
     \label{Figure 1}
\end{center}
\end{figure}

\begin{figure}[htb]
     \begin{center}
     \begin{tabular}{c}
     \includegraphics[width=10cm]{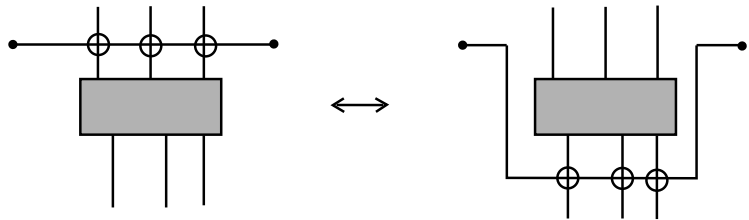}
     \end{tabular}
     \caption{\bf Detour Move}
     \label{Figure 2}
\end{center}
\end{figure}

\begin{figure}[htb]
     \begin{center}
     \begin{tabular}{c}
     \includegraphics[width=10cm]{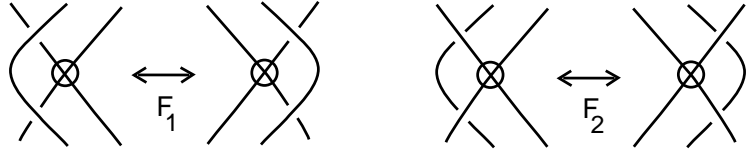}
     \end{tabular}
     \caption{\bf Forbidden Moves}
     \label{Figure 3}
\end{center}
\end{figure}

Another way to understand virtual diagrams is to regard them as representatives for oriented Gauss codes \cite{GPV}, \cite{VKT,SVKT} 
(Gauss diagrams). Such codes do not always have planar realizations. An attempt to embed such a code in the plane
leads to the production of the virtual crossings. The detour move makes the particular choice of virtual crossings 
irrelevant. {\it Virtual isotopy is the same as the equivalence relation generated on the collection
of oriented Gauss codes by abstract Reidemeister moves on these codes.}  
\bigbreak

Figure $3$ illustrates the two {\it forbidden moves}. Neither of these follows from Reidmeister moves plus detour move, and 
indeed it is not hard to construct examples of virtual knots that are non-trivial, but will become unknotted on the application of 
one or both of the forbidden moves. The forbidden moves change the structure of the Gauss code and, if desired, must be 
considered separately from the virtual knot theory proper. 
\bigbreak

\section{Interpretation of Virtuals Links as Stable Classes of Links in  Thickened Surfaces}
There is a useful topological interpretation \cite{VKT,DVK} for this virtual theory in terms of embeddings of links
in thickened surfaces.  Regard each 
virtual crossing as a shorthand for a detour of one of the arcs in the crossing through a 1-handle
that has been attached to the 2-sphere of the original diagram.  
By interpreting each virtual crossing in this way, we
obtain an embedding of a collection of circles into a thickened surface  $S_{g} \times R$ where $g$ is the 
number of virtual crossings in the original diagram $L$, $S_{g}$ is a compact oriented surface of genus $g$
and $R$ denotes the real line.  We say that two such surface embeddings are
{\em stably equivalent} if one can be obtained from another by isotopy in the thickened surfaces, 
homeomorphisms of the surfaces and the addition or subtraction of empty handles (i.e. the knot does not go through the handle).

\begin{figure}
     \begin{center}
     \begin{tabular}{c}
     \includegraphics[width=10cm]{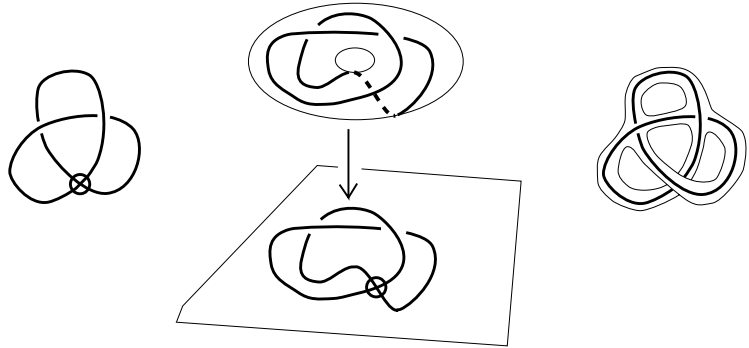}
     \end{tabular}
     \caption{\bf Surfaces and Virtuals}
     \label{Figure 4}
\end{center}
\end{figure}

\begin{figure}
     \begin{center}
     \begin{tabular}{c}
     \includegraphics[width=10cm]{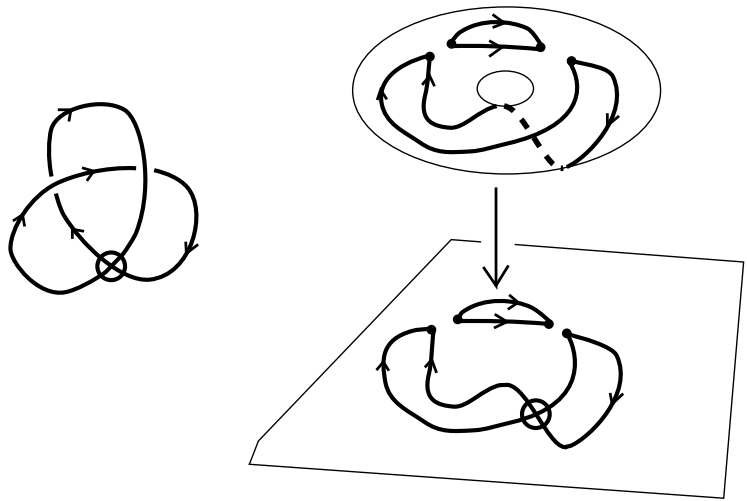}
     \end{tabular}
     \caption{\bf Surfaces and Virtual States}
     \label{Figure 5}
\end{center}
\end{figure}

\noindent We have the
\smallbreak
\noindent
{\bf Theorem 1 \cite{VKT,DKT,DVK,Carter}.} {\em Two virtual link diagrams are isotopic if and only if their corresponding 
surface embeddings are stably equivalent.}  
\smallbreak
\noindent
\bigbreak  

\noindent In Figure 4 we illustrate some points about this association of virtual diagrams and knot and link diagrams on surfaces.
Note the projection of the knot diagram on the torus to a diagram in the plane (in the center of the figure) has a virtual crossing in the 
planar diagram where two arcs that do not form a crossing in the thickened surface project to the same point in the plane. In this way, virtual 
crossings can be regarded as artifacts of projection. The same figure shows a virtual diagram on the left and an ``abstract knot diagram" \cite{Kamada3,Carter} on the right.
The abstract knot diagram is a realization of the knot on the left in a thickened surface with boundary and it is obtained by making a neighborhood of the 
virtual diagram that resolves the virtual crossing into arcs that travel on separate bands. The virtual crossing appears as an artifact of the
projection of this surface to the plane. The reader will find more information about this correspondence \cite{VKT,DKT} in other papers by the author and in
the literature of virtual knot theory.
\bigbreak
 
\section{Flat Virtual Knots and Links}
Every classical knot or link diagram can be regarded as a $4$-regular plane graph with extra structure at the 
nodes. This extra structure is usually indicated by the over and under crossing conventions that give
instructions for constructing an embedding of the link in three dimensional space from the diagram.  If we take the flat diagram
without this extra structure then the diagram is the shadow of some link in three dimensional space, but the weaving of that link is not 
specified. It is well known that if one is allowed to apply the Reidemeister moves to such a shadow (without regard to the types
of crossing since they are not specified) then the shadow can be reduced to a disjoint union of circles. This reduction is 
no longer true for virtual links. More precisely, let a {\em flat virtual diagram} be a diagram with {\it virtual crossings} as we have
described them and {\em flat crossings} consisting in undecorated nodes of the $4$-regular plane graph. Two flat virtual diagrams are {\em equivalent} if
there is a  sequence of generalized flat Reidemeister moves (as illustrated in Figure 1) taking one to the other. A generalized
flat Reidemeister move is any move as shown in Figure 1 where one ignores the over or under crossing structure.
Note that in studying flat virtuals the rules for changing virtual crossings among themselves and the rules for changing
flat crossings among themselves are identical. Detour moves as in Figure 1C are available for virtual crossings
with respect to flat crossings and {\it not} the other way around. The analogs of the forbidden moves of Figure 3 remain forbidden 
when the classical crossings are replaced by flat crossings. 
\bigbreak

The theory of flat virtual knots and links is identical to the theory of all oriented Gauss codes (without over or under information)
modulo the flat Reidemeister moves. Virtual crossings are an artifact of the realization of the flat diagram in the plane.
In Turaev's work \cite{VST} flat virtual knots and links are called {\it virtual strings}. See also recent papers of Roger Fenn \cite{Fenn1,Fenn2} for other
points of view about flat virtual knots and links.
\bigbreak

We shall say that a virtual diagram {\em overlies} a flat diagram if the virtual diagram is obtained from the flat diagram by
choosing a crossing type for each flat crossing in the virtual diagram. To each virtual diagram $K$ there is an associated 
flat diagram $F(K)$, obtained by forgetting the extra structure at the classical crossings in $K.$ Note that if $K$ and $K'$
are isotopic as virtual diagrams, then $F(K)$ and $F(K')$ are isotopic as flat virtual diagrams. Thus, if we can
show that $F(K)$ is not reducible to a disjoint union of circles, then it will follow that $K$ is a non-trivial virtual link.
The flat virtual diagrams present a challenge for the
construction of new invariants. They are fundamental to the study of virtual knots.  A virtual knot is necessarily non-trivial if its flat projection
is a non-trivial flat virtual diagram.  We wish to be able to determine when a given virtual  link is isotopic to a classcal
link. The reducibility or irreducibility of the underlying flat diagram is the first  obstruction to such an equivalence.
\bigbreak 

\noindent {\bf Definition.}
A {\it virtual graph} is a flat virtual diagram where the classical flat crossings are not subjected to the flat Reidemeister moves.
Thus a virtual graph is a $4-regular$ graph that is represented in the plane via a choice of cyclic orders at its nodes. The virtual crossings
are artifacts of this choice of placement in the plane, and we allow detour moves for consecutive sequences of virtual crossings just as in the
virtual knot theory. Two virtual graphs are {\it isotopic} if there is a combination of planar graph isotopies and detour moves that connect them.
The theory of virtual graphs is equivalent to the theory of $4$-regular graphs on oriented surfaces, taken up to empty handle stabilization, in 
direct analogy to our description of the theory of virtual links and flat virtual links.
\bigbreak

\section{Long Virtual Knots, Long Flats, Parity and the Odd Writhe}
A long knot or link is a $1-1$ tangle. It is a tangle with one input end and one output end. In between one has, in the diagram, 
any knotting or linking, virtual or classical. Classical long knots (one component) carry essentially the same topological information
as their closures. In particular, a classical long knot is knotted if and only if its closure is knotted. This statement is false for virtual
knots. An example of the phenomenon is shown in Figure 5.
\bigbreak

The long knots $L$ and $L'$ shown in Figure 5 are non-trivial in the virtual category. Their closures, obtained by attaching the ends together are
unknotted virtuals. Concomittantly, there can be a multiplicity of long knots associated to a given virtual knot diagram, obtained by cutting an
arc from the diagram and creating a $1-1$ tangle. {\it It is a fundamental problem to determine the kernel of the closure mapping from long
virtual knots to virtual knots.} In Figure 5, the long knots $L$ and $L'$ are trivial as welded long knots (where the first forbidden move is allowed). The
obstruction to untying them as virtual long knots comes from the first forbidden move. The matter of proving that $L$ and $L'$ are non-trivial distinct long
knots is difficult by direct attack. There is a fundamental relationship betweem long flat virtual knots and long virtual knots that can be used to
see it.
\bigbreak

Let $LFK$ denote the set of long flat virtual knots and let $LVK$ denote the set of long virtual knots.
We define $$A:LFK \longrightarrow LVK$$ by letting $A(S)$ be the ascending long virtual knot diagram associated with 
the long flat virtual diagram $S.$ That is, $A(S)$ is obtained from $S$ by traversing $S$ from its left end to its right end and creating
a crossing at each flat crossing so that {\it one passes under each crossing before passing over that crossing}. Virtual crossings are not changed by this
construction. The idea of using the ascending diagram to define invariants of long flat virtuals is exploited in \cite{SW}. The following result is new.
\bigbreak

\noindent {\bf Long Flat Embedding Theorem 2.} The mapping $A:LFK \longrightarrow LVK$ is well-defined on the corresponding isotopy classes of diagrams, 
and it is injective. Hence {\it long flat virtual knots embed in the class of long virtual knots.}
\bigbreak

\noindent {\bf Proof.} It is easy to see that if two flat long diagrams $S,T$ are virtually isotopic then $A(S)$ and $A(T)$ are isotopic 
long virtual knots. We leave the verification to the reader. To see the injectivity, define 
$$Flat:LVK \longrightarrow LFK$$ by letting $Flat(K)$ be the long flat diagram obtained from $K$ by flattening all the classical crossings in $K.$
By definition $Flat(A(S)) = S.$ Note that the map $Flat$ takes isotopic long virtual knots to isotopic flat virtual knots since
any move using the under and over crossings is a legitimate move when this distinction is forgotten. This proves injectivity and
completes the proof of the Theorem. //
\bigbreak

The long flat theorem is easy to prove, and it has many good consequences. First of all, let $Inv(K)$ denote any invariant
of long virtual knots. ($Inv(K)$ can denote a polynomial invariant, a number, a group, quandle, biquandle or other invariant structure.) Then we can define
$$Inv(S)$$ for any long flat knot $S$ by the formula $$Inv(S) = Inv(A(S)),$$ and this definition yields an invariant of long flat knots.
\bigbreak

\begin{figure}
     \begin{center}
     \begin{tabular}{c}
     \includegraphics[width=8cm]{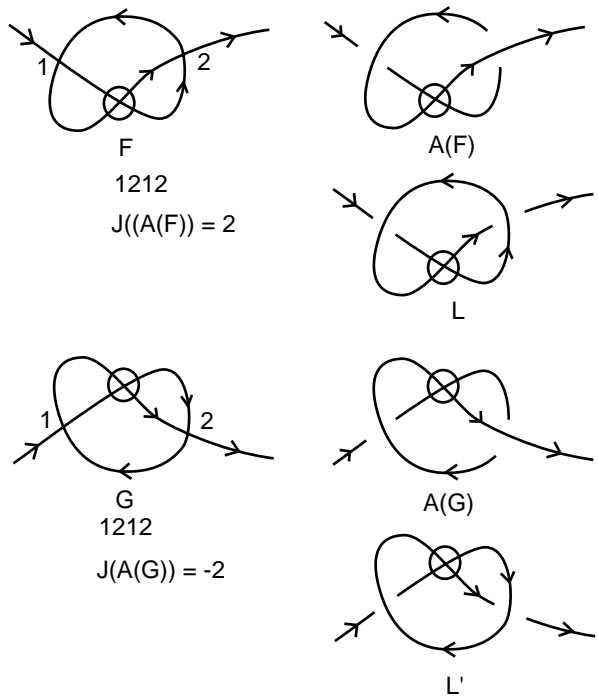}
     \end{tabular}
     \caption{\bf Ascending Map}
     \label{Figure 6}
\end{center}
\end{figure}

\noindent {\bf Using the Odd Writhe J(K).} View Figure 5. We show the long flat $F$ and its image under the
ascending map, $A(F)$ are non-trivial. In fact, $A(F)$ is non-trivial and non-classical. One computes that $J(A(F))$ is
non-zero where $J(K)$ denotes the {\it odd writhe} of $K.$ The odd writhe \cite{SL} is the sum of the signs of the odd crossings. A crossing is {\it odd} if 
it flanks an odd number of symbols in the Gauss code of the diagram. Classical diagrams have zero odd writhe. 
In this case, the flat Gauss code for $F$ is $1212$ with both  crossings odd. Thus we see from the figure that $J(A(F)) = 2.$
Thus $A(F)$ is non-trivial, non-classical and inequivalent to its mirror image.
Once we check that $A(F)$ is non-trivial, we know that the flat knot $F$ is non-trivial, and from this we conclude that the long virtual
$L$ is also non-trivial. Note that $J(L) = 0,$ so we cannot draw this last conclusion directly from $J(L).$ This same figure illustrates a long flat $G$
that is obtained by reflecting $F$ in a horizontal line. Then, as the reader can calculate from this figure, $J(A(G)) = -2.$ Thus $F$ and $G$ 
are distinct non-trivial long flats. We conclude from these arguments that the long virtual knots $L$ and $L'$ in Figure 5 are both non-trivial, and that
$L$ is not virtually isotopic to $L'$ (since such an isotopy would give an isotopy of $F$ with $G$ by the flattening map).
\bigbreak

\section{Review of the Bracket Polynomial for Virtual Knots}

In this section we recall how the bracket state summation model \cite{K} for the Jones polynomial \cite{Jones,WITT} is defined for virtual knots
and links.  In the next section we give an extension of this model using orientation structures on the states of the bracket expansion.
The extension is also an invariant of flat virtual links.
\bigbreak

We call a diagram in the plane 
{\em purely virtual} if the only crossings in the diagram are virtual crossings. Each purely virtual diagram is equivalent by the
virtual moves to a disjoint collection of circles in the plane.
\bigbreak

A state $S$ of a link diagram $K$ is obtained by
choosing a smoothing for each crossing in the diagram and labelling that smoothing with either $A$ or $A^{-1}$
according to the convention that a counterclockwise rotation of the overcrossing line sweeps two 
regions labelled $A$, and that a smoothing that connects the $A$ regions is labelled by the letter $A$. Then, given
a state $S$, one has the evaluation $<K|S>$ equal to the product of the labels at the smoothings, and one has the 
evaluation $||S||$ equal to the number of loops in the state (the smoothings produce purely virtual diagrams).  One then has
the formula
$$<K> = \Sigma_{S}<K|S>d^{||S||-1}$$
where the summation runs over the states $S$ of the diagram $K$, and $d = -A^{2} - A^{-2}.$
This state summation is invariant under all classical and virtual moves except the first Reidemeister move.
The bracket polynomial is normalized to an
invariant $f_{K}(A)$ of all the moves by the formula  $f_{K}(A) = (-A^{3})^{-w(K)}<K>$ where $w(K)$ is the
writhe of the (now) oriented diagram $K$. The writhe is the sum of the orientation signs ($\pm 1)$ of the 
crossings of the diagram. The Jones polynomial, $V_{K}(t)$ is given in terms of this model by the formula
$$V_{K}(t) = f_{K}(t^{-1/4}).$$
\noindent This definition is a direct generalization to the virtual category of the  
state sum model for the original Jones polynomial. It is straightforward to verify the invariances stated above.
In this way one has the Jones polynomial for virtual knots and links.
\bigbreak

\noindent We have \cite{DVK} the  
\smallbreak
\noindent
{\bf Theorem 3.} {\em To each non-trivial
classical knot diagram of one component $K$ there is a corresponding  non-trivial virtual knot diagram $Virt(K)$ with unit
Jones polynomial.} 
\bigbreak

\noindent {\bf Proof Sketch.} This Theorem is a key ingredient in the problems involving virtual knots. Here is a sketch of its proof.
The proof uses two invariants of classical knots and links that generalize to arbitrary virtual knots and links.
These invariants are the {\em Jones polynomial} and the {\em involutory quandle} denoted by the notation
$IQ(K)$ for a knot or link $K.$ 
\bigbreak

Given a
crossing $i$ in a link diagram, we define $s(i)$ to be the result of {\em switching} that crossing so that the undercrossing arc
becomes an overcrossing arc and vice versa. We define the {\em virtualization}
$v(i)$ of the crossing by the local replacement indicated in Figure 7. In this figure we illustrate how, in virtualization, 
the  original crossing is replaced by a crossing that is flanked by two virtual crossings. When we smooth the two
virtual crossings in the virtualization we obtain the original knot or link diagram with the crossing switched.
\bigbreak

Suppose that $K$ is a (virtual or classical) diagram with a classical crossing labeled $i.$  Let $K^{v(i)}$ be the diagram
obtained from $K$ by virtualizing the crossing $i$ while leaving the rest of the diagram just as before. Let $K^{s(i)}$ be
the diagram obtained from $K$ by switching the crossing $i$ while leaving the rest of the diagram just as before. Then it
follows directly from the expansion formula for the bracket polynomial that $$V_{K^{s(i)}}(t) = V_{K^{v(i)}}(t).$$ 
\noindent As far as the Jones
polynomial is concerned, switching a crossing and virtualizing a crossing look the same. We can start with a classical knot diagram $K$
and choose a subset $S$ of crossings such that the diagram is unknotted when these crossings are switched. Letting $Virt(K)$ denote the virtual knot
diagram obtained by virtualizing each crossing in $S$, it follows that the Jones polynomial of $Virt(K)$ is equal to unity, the Jones polynomial
of the unknot. Nevertheless, if the original knot $K$ is knotted, then the virtual knot $Virt(K)$ will be non-trivial. We outline the argument for this
fact below.
\bigbreak 
 
The involutory quandle \cite{KNOTS} is an algebraic invariant
equivalent to the fundamental group of the double branched cover of a knot or link in the classical case. In this algebraic
system one associates a generator of the algebra $IQ(K)$ to each arc of the diagram $K$ and there is a relation of the form
$c = ab$ at each crossing, where $ab$ denotes the (non-associative) algebra product of $a$ and $b$ in $IQ(K).$ See Figure 8.
In this figure we have illustrated the fact that  $$IQ(K^{v(i)}) = IQ(K).$$
As far as the involutory quandle is concerned, the original crossing and the virtualized crossing look the same.
\bigbreak

If a classical knot is actually knotted, then its involutory quandle is non-trivial \cite{W}. Hence if we start
with a non-trivial classical knot and virtualize any subset of its crossings we obtain a virtual knot that is still 
non-trivial. There is a subset $A$ of the crossings of a classical knot $K$ such that the knot $SK$ obtained by
switching these crossings is an unknot.  Let $Virt(K)$ denote the virtual diagram obtained from $A$ by virtualizing
the crossings in the subset $A.$  By the above discussion, the Jones polynomial of $Virt(K)$ is the same
as the Jones polynomial of $SK$, and this is $1$ since $SK$ is unknotted. On the other hand, the $IQ$ of $Virt(K)$ is the 
same as the $IQ$ of $K$, and hence if $K$ is knotted, then so is $Virt(K).$   We have shown that $Virt(K)$ is a non-trivial
virtual knot with unit Jones polynomial.  This completes the proof of the Theorem. //
\bigbreak

\begin{figure}
     \begin{center}
     \begin{tabular}{c}
     \includegraphics[width=8cm]{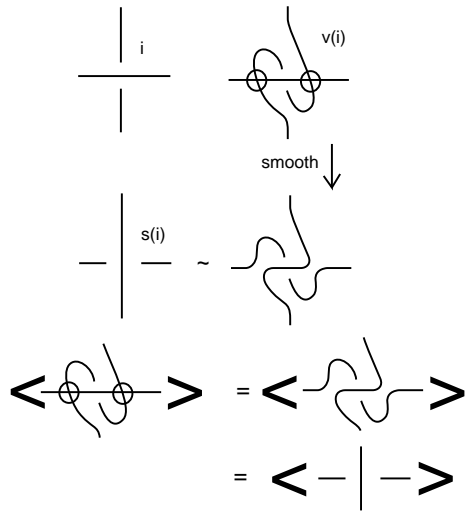}
     \end{tabular}
     \caption{\bf Switch and Virtualize}
     \label{Figure 7}
\end{center}
\end{figure}

\begin{figure}
     \begin{center}
     \begin{tabular}{c}
     \includegraphics[width=8cm]{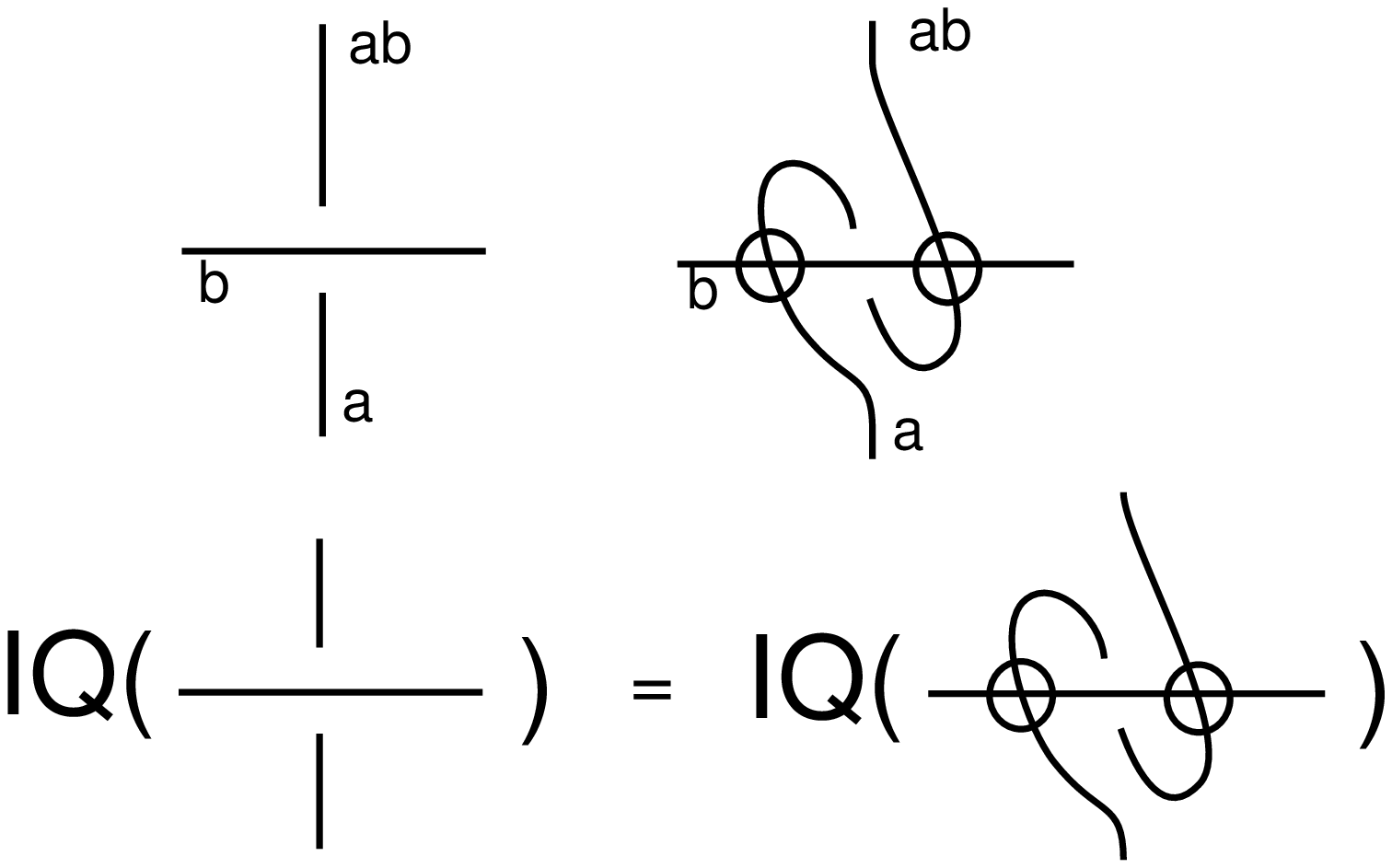}
     \end{tabular}
     \caption{\bf IQ(Virt)}
     \label{Figure 8}
\end{center}
\end{figure}

See Figure 32 for an example of a virtualized trefoil, the simplest example of a non-trivial virtual knot with unit Jones polynomial.
More work is needed to prove that the virtual knot $T$ in Figure 32 is not classical. In the next section we will give a proof of this fact
by using an extension of the bracket polynomial.
\bigbreak

It is an open problem whether there are classical knots (actually knotted) having unit Jones polynomial. (There are linked links whose linkedness is unseen
\cite{MT,EKT} by the Jones polynomial.) If there exists a classical knot with unit Jones polynomial, then one of the knots $Virt(K)$ produced by this Theorem
may be isotopic to  a classical knot.  Such examples are guaranteed to be non-trivial, but they are usually also not classical. 
We do not know at this writing whether all such virtualizations of non-trivial classical knots, yielding virtual knots with unit Jones polynomial, are 
non classical. It is
an intricate task to verify that specific examples of 
$Virt(K)$ are not classical. This has led to an investigation of new invariants for virtual knots. In this way the search for classical knots with unit Jones
polynomial expands to exploration of the structure of the infinite collection of virtual knots with unit Jones polynomial.
\bigbreak

\begin{figure}
     \begin{center}
     \begin{tabular}{c}
     \includegraphics[width=3cm]{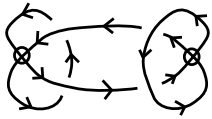}
     \end{tabular}
     \caption{\bf Kishino Diagram}
     \label{Figure 9}
\end{center}
\end{figure}

In Figure 9 we show the {\it Kishino diagram} $K.$
This diagram has unit Jones polynomial and its fundamental group is infinite cyclic. The Kishino diagram was discovered by Kishino in \cite{KIS}.
Many other invariants of virtual knots 
fail to detect the Kishino knot. Thus it has been a test case for examining new invariants. Heather Dye and the author \cite{MinSurf}
have used the bracket polynomial defined for knots and links in a thickened surface (the state curves
are taken as isotopy classes of curves in the surface) to prove the non-triviality and non-classicality of the Kishino diagram.
In fact, we have used this technique to show that knots with unit Jones polynomial obtained by a single virtualization are non-classical.
See the problem list by Fenn, Kauffman and Manturov \cite{VP} for other problems and proofs related to the Kishino diagram.
In the next section we describe a new extension of the bracket polynomial that can be used to discriminate the Kishino diagram, and, in fact,
shows that its corresponding flat virtual knot is non-trivial.
\bigbreak

\section{The Parity Bracket Polynomial}
In this section we introduce the Parity Bracket Polynomial of Vassily Manturov \cite{MP}.
This is a generalization of the bracket polynomial to virtual knots and links that uses the parity of the crossings. We define a {\em Parity State} of a virtual diagram $K$ to be a labeled virtual
graph obtained from $K$ as follows: For each odd crossing in $K$ replace the crossing by a graphical node. For each even crossing in $K$ replace the crossing by one of its two possible smoothings, and label the smoothing site by $A$ or $A^{-1}$ in the usual way. Then we define the parity bracket by the 
state expansion formula
$$\langle K \rangle _{P} = \sum_{S}A^{n(S)}[S]$$
where $n(S)$ denotes the number of $A$-smoothings minus the number of $A^{-1}$ smoothings and 
$[S]$ denotes a combinatorial evaluation of the state defined as follows: First reduce the state
by Reidemeister two moves on nodes as shown in Figure 10. The graphs are taken up to virtual equivalence (planar isotopy plus detour moves on the virtual crossings). Then regard the reduced state
as a disjoint union of standard state loops (without nodes) and graphs that irreducibly contain nodes.
With this we write $$[S] = (- A^{2} - A^{-2})^{l(S)} [G(S)]$$ where $l(S)$ is the number of standard loops in the reduction of the state $S$ and $[G(S)]$ is the disjoint union of reduced graphs that contain nodes.
In this way, we obtain a sum of Laurent polynomials in $A$ multiplying reduced graphs as the Manturov Parity Bracket. It is not hard to see that this bracket is invariant under regular isotopy and detour moves
and that it behaves just like the usual bracket under the first Reidemeister move. However, the use of
parity to make this bracket expand to graphical states gives it considerable extra power in some situations. For example, consider the Kishino diagram in Figure 9. We see that all the classical crossings in this knot are odd. Thus the parity bracket is just the graph obtained by putting nodes at each of these crossings. The resulting graph does not reduce under the graphical Reidemeister two moves, and so we conclude that the Kishino knot is non-trivial and non-classical. Since we can apply the parity
bracket to a flat knot by taking $A = -1$, we see that this method shows that the Kishino flat is non-tirival. 
\bigbreak

\begin{figure}
     \begin{center}
     \begin{tabular}{c}
     \includegraphics[width=8cm]{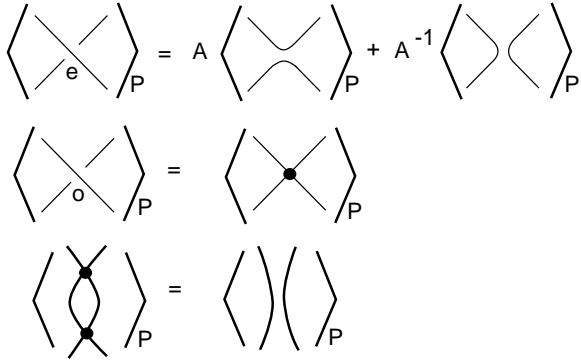}
     \end{tabular}
     \caption{\bf Parity Bracket Expansion}
     \label{Figure 10}
\end{center}
\end{figure}

\begin{figure}
     \begin{center}
     \begin{tabular}{c}
     \includegraphics[width=4cm]{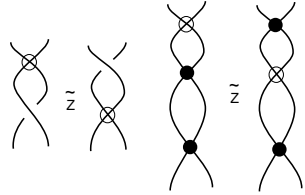}
     \end{tabular}
     \caption{\bf Z-Move and Graphical Z-Move}
     \label{Figure 11}
\end{center}
\end{figure}

\begin{figure}
     \begin{center}
     \begin{tabular}{c}
     \includegraphics[width=6cm]{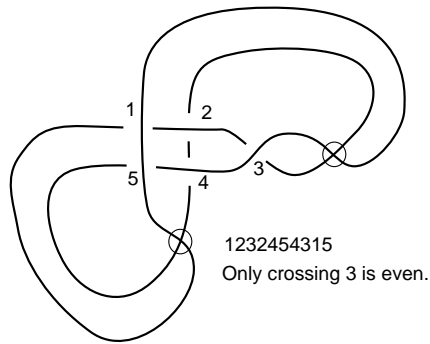}
     \end{tabular}
     \caption{\bf A Knot KS With Unit Jones Polynomial}
     \label{Figure 12}
\end{center}
\end{figure}

\begin{figure}
     \begin{center}
     \begin{tabular}{c}
     \includegraphics[width=6cm]{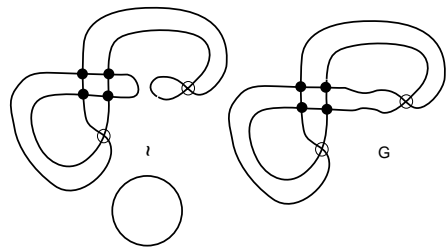}
     \end{tabular}
     \caption{\bf Parity Bracket States for the Knot KS}
     \label{Figure 13}
\end{center}
\end{figure}

In Figure 11 we illustrate the {\it Z-move}  and the {\it graphical Z-move}. Two virtual knots or links
that are related by a Z-move have the same standard bracket polynomial. This follows directly from our discussion in the previous section. We would like to analyze the structure of Z-moves using the parity
bracket. In order to do this we need a version of the parity bracket that is invariant under the Z-move.
In order to accomplish this, we need to add a corresponding Z-move in the graphical reduction process for the parity bracket. This extra graphical reduction is indicated in Figure 11 where we show a graphical Z-move. The reader will note that graphs that are irreducible without the graphical Z-move can become reducible if we allow graphical Z-moves in the reduction process. For example, the graph associated with the Kishino knot is reducible under graphical Z-moves. However, there are examples of 
graphs that are not reducible under graphical Z-moves and Reidemister two moves. An example of such a graph occurs in the parity bracket of the knot $KS$ shown in Figure 12. This knot has one even classical crossing and four odd crossings. One smoothing of the even crossing yields a state that reduces to a loop with no graphical nodes, while the other smoothing yields a state that is irreducible even when the Z-move is allowed. The upshot is that this knot KS is not Z-equivalent to any classical knot. Since one can verify that $KS$ has unit Jones polynomial, this example is a counterexample to a conjecture of Fenn, Kauffman and Maturov \cite{VP} that suggested that a knot with unit Jones polynomial should be Z-equivalent to a classical knot.
\bigbreak

Parity is clearly an important theme in virtual knot theory and will figure in many future investigations of this subject. The type of construction that we have indicated for the bracket polynomial in this section
can be varied and applied to other invariants. Furthermore the notion of describing a parity for crossings
in a diagram is also susceptible to generalization. For more on this theme the reader should consult
\cite{MP1,PT} and \cite{SL} for our original use of parity for another variant of the bracket polynomial.
\bigbreak

\noindent {\bf Remark on Free Knots}
Manturov (See \cite{MP}) has defined the domain of {\it free knots}. A free knot is a Gauss diagram
(or Gauss code) without any orientations or signs, taken up to abstract Reidemeister moves. The reader 
will see easily that free knots are the same as flat virtual knots modulo the flat $Z$-move. Furthermore,
the parity bracket polynomial evaluated at $A = 1$ or $A = -1$ is an invariant of free knots. By using it
on examples where all the crossings are odd, one obtains infinitely many examples of non-trivial free knots. This is just the beginning of an investigation into free knots that will surely lead to deep relationsips between combinatorics and knot theory. Any free knot that is shown to be non-trivial has
a number of non-trivial virtual knots overlying it. The free knots deserve to be studied for their own sake.
A first example of new invariants of free knots that go beyond Gaussian parity can be found in the paper
\cite{KMKup}.
\bigbreak

\section{The Arrow Polynomial for Virtual and Flat Virtual Knots and Links}
This section describes an invariant for oriented virtual knots and links,  and for flat oriented virtual knots and links that we call the {\em arrow polynomial} \cite{DyeKauff,ExtBr}. This invariant is considerably 
stronger than the Jones polynomial for virtual knots and links, and it is a very natural extension of
the Jones polynomial, using the oriented diagram structure of the state summation.
The construction of the arrow polynomial invariant begins with the
oriented state summation of the bracket polynomial. This means that each local smoothing is either an oriented smoothing or a {\it disoriented
smoothing} as illustrated in Figures 14 and 15.
In \cite{DyeKauff} we  show how the arrow polynomial  can be used to estimate virtual crossing numbers.
In \cite{ExtBr} we  discuss the arrow polynomial and also a generalization, the {\em extended bracket polynomial}. The extended bracket polynomial  is harder to compute than the arrow polynomial, and we
shall not discuss its properties in this introduction. We recommend that the interested reader consult our
 paper \cite{ExtBr} on this subject. We describe a new extension of the arrow polynomial to long virtual knots at the end of the present section.
\bigbreak

\bigbreak

\begin{figure}
     \begin{center}
     \begin{tabular}{c}
     \includegraphics[width=6cm]{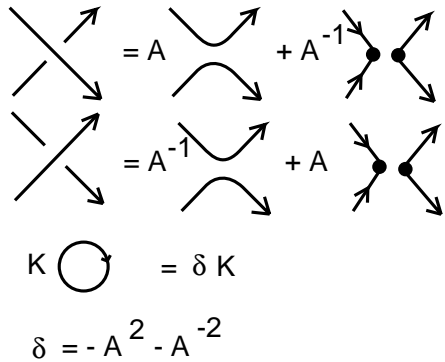}
     \end{tabular}
     \caption{\bf Oriented Bracket Expansion}
     \label{Figure 14}
\end{center}
\end{figure}

\begin{figure}
     \begin{center}
     \begin{tabular}{c}
     \includegraphics[width=6cm]{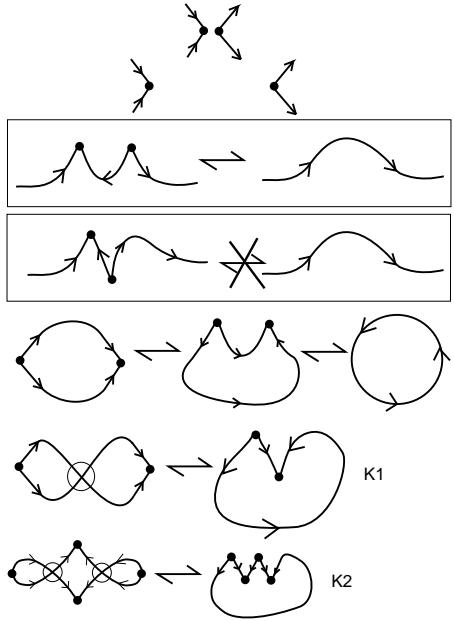}
     \end{tabular}
     \caption{\bf Reduction Relation for the Arrow Polynomial.}
     \label{Figure 15}
\end{center}
\end{figure}

\begin{figure}
     \begin{center}
     \begin{tabular}{c}
     \includegraphics[width=8cm]{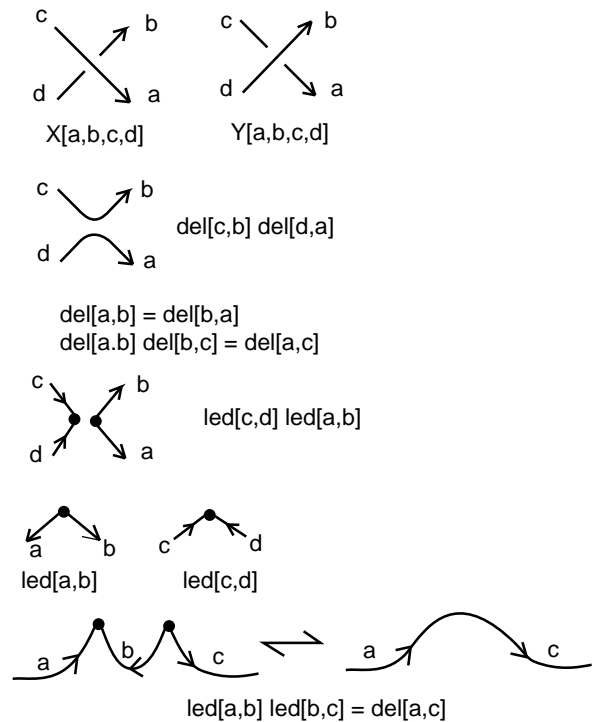}
     \end{tabular}
     \caption{\bf Formal Kronecker Deltas.}
     \label{Figure 16}
\end{center}
\end{figure}

\begin{figure}
     \begin{center}
     \begin{tabular}{c}
     \includegraphics[width=8cm]{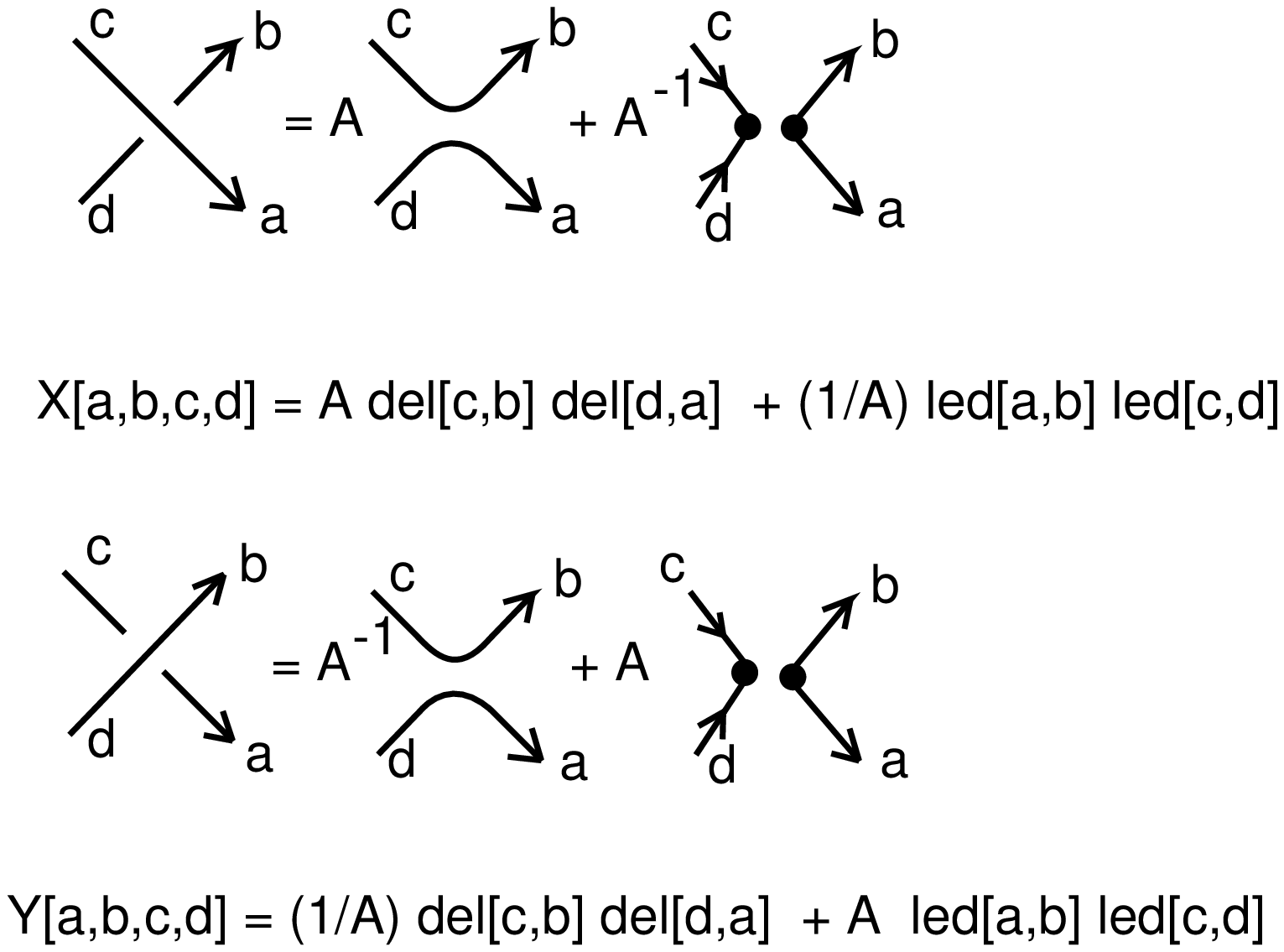}
     \end{tabular}
     \caption{\bf Simple arrow polynomial Expansion via Formal Kronecker Deltas.}
     \label{Figure 17}
\end{center}
\end{figure}

In Figure 14 we illustrate the oriented bracket expansion for both positive and negative crossings in a  link diagram. An oriented crossing can
be smoothed in the oriented fashion or the disoriented fashion as shown in Figure 14. We refer to these smoothings as {\it oriented} and {\it disoriented}
smoothings. To each smoothing we make an associated configuration that will be part of the arrow polynomial state summation.  The configuration associated to a state with oriented and  disoriented smoothings is obtained by applying the reduction rules described below.  See Figures 15 and 16. The
arrow polynomial state summation is defined by the formula:
$$\langle \langle K \rangle \rangle = {\cal A}[K] = \Sigma_{S}\langle K|S \rangle d^{||S||-1} [S]$$
where $S$ runs over the oriented bracket states of the diagram, $\langle K|S \rangle $ is the usual product of vertex weights as in the  standard bracket polynomial, and $[S]$ is a product of extra variables $K_{1}, K_{2}, \cdots$ associated with the state $S.$ These variables are explained below.
Note that we use the designation $\langle \langle K \rangle \rangle$ to indicate the arrow polynomial 
in this paper and in the figures in this paper. In \cite{ExtBr} we have used this notation to refer to the 
extended bracket generalization of the arrow polynomial. Since we do not use the extended bracket in this paper, it is convenient to use the double bracket notation here for the arrow polynomial.
\bigbreak

Due to the oriented state expansion, the loops in the resulting states have extra combinatorial structure in the form of paired {\it cusps} as shown in Figure 14. Each disoriented smoothing gives rise to 
a cusp pair where each cusp has either two oriented lines going into the cusp or two oriented lines leaving the cusp. We reduce this structure according to a set of rules that yields invariance of the state summation under the Reideimeister moves. The basic conventions for this 
simplification are shown in Figures 15 and 16.  Each cusp is denoted by an angle with arrows either both entering the vertex or both leaving the vertex. Furthermore, the angle locally divides the 
plane into two parts: One part is the span of an acute angle; the other part is the span of an obtuse angle. We refer to the span of the acute angle as
the {\it inside} of the cusp.  
\bigbreak

\noindent {\bf Remark on State Reduction.}  View Figures 15 and 16.  The basic reduction move in these figures corresponds to the elimination of two consecutive cusps on a single loop.
Note that we allow cancellation of consecutive cusps along a loop where the cusps both 
point to the same local region of the two local regions delineated by the loop. We do not allow cancellation of a  ``zig-zag" where two consecutive cusps point to opposite (local) sides of the loop.      We shall note that in a classical knot or link diagram, all state loops reduce to
loops that are free from cusps.  
\bigbreak

Figure 15 illustates the basic reduction rule for the {\it arrow polynomial.}  The reduction rule allows the cancellation of two adjacent cusps when they have {\it insides on the same
side} of the segment that connects them. When the insides of the cusps are on opposite sides of the connecting segment, then no cancellation is allowed. Each state circle is seen as a circle graph with extra nodes corresponding to the cusps.
All graphs are taken up to virtual equivalence, as explained earlier in this paper. Figure 15 illustrates the simplification of two circle graphs. In one case
the graph reduces to a circle with no vertices. In the other case there is no further cancellation, but the graph is equivalent to one without a virtual crossing.
The state expansion for ${\cal A}[K] = <<K>>$ is exactly as shown in Figure 14, but we use the reduction rule of Figure 15 so that each state is a disjoint union of reduced circle 
graphs. Since such graphs are planar, each is equivalent to an embedded graph (no virtual crossings) 
via the detour move, and the reduced forms of such graphs have $2n$ vertices that 
alternate in type around the circle so that $n$ are pointing inward and $n$ are pointing outward. The circle with no vertices is evaluated as $d = -A^2 - A^{-2}$ as
is usual for these expansions, and the circle is removed from the graphical expansion. We let $K_{n}$ denote the circle graph with $2n$ alternating vertex types as shown in
Figure 15 for $n=1$ and $n=2.$
Each circle graph contributes $d = -A^2 - A^{-2}$ to the state sum and the graphs $K_{n}$ for $n \ge 1$ remain in the graphical expansion. Each $K_{n}$ is an extra variable in the polynomial. Thus a product of the $K_{n}$'s
corresponds to a state that is a disjoint union of copies of these circle graphs. By evaluating each circle graph as $d = -A^2 - A^{-2}$ (as well as taking its arrow variable $K_{n}$) we guarantee
that the resulting polynomial will reduce to the original bracket polynomial when each of the new variables $K_{n}$ is set equal to unity. Note that we 
continue to use the caveat that an isolated circle or circle graph (i.e. a state consisting in a single circle or single circle graph) is assigned a loop value
of unity in the state sum. This assures that ${\cal A}[K]$ is normalized so that the unknot receives the value one.
\bigbreak

We have the following Proposition, showing that the phenomenon of cusped states and extra variables
$K_{n}$ only occurs for virtual knots and links.
\bigbreak

\noindent {\bf Proposition 2.} In a classical knot or link diagram, all state loops reduce to
loops that are free from cusps. 
\smallbreak

\noindent {\bf Proof.} The result follows from the Jordan Curve Theorem. Each state loop for a classical knot is a Jordan curve on the surface of the 
$2$-sphere, dividing the surface of the sphere into an inside and an outside. If a state loop has a non-trivial arrow reduction, then there will be non-empty collections of inward-pointing cusps and a collection of outward-pointing cusps
in the reduced loop. Each cusp must be paired with another cusp in the state (since cusps are originally paired). No two inward-pointing cusps on a loop can be paired with one another, since they will have incompatible
orientations (similarly for the outward-pointing cusps). Therefore, for a given non-trivial reduced loop $\lambda$ in a classical diagram there must be another non-trivially
reduced loop $\lambda'$ inside  $\lambda$ and another such loop $\lambda''$ outside $\lambda$ to handle the necessary pairings. This means that the given state
of the diagram would have infinitely many loops. Since we work with knot and link diagrams with finitely many crossings, this is not possible. Hence there are no
non-trivially reduced loops in the states of a classical diagram. This completes the proof. //
\bigbreak

\noindent We now have the following state summation for the arrow polynomial $${\cal A}[K] = \Sigma_{S}<K|S>d^{||S||-1} {\cal V}[S]$$
where $S$ runs over the oriented bracket states of the diagram, $<K|S>$ is the usual product of vertex weights as in the 
standard bracket polynomial, $||S||$ is the number of circle graphs in the state $S$, and ${\cal V}[S]$ is a product of the variables $K_{n}$ associated
with the non-trivial circle graphs in the state $S.$ Note that each circle graph (trivial or not) contributes to the power of $d$ in the state summation,
but only non-trivial circle graphs contribute to ${\cal V}[S].$ 
\bigbreak

\noindent {\bf Theorem 6.} With the above conventions, the arrow polynomial ${\cal A}[K]$ is a polynomial in $A, A^{-1}$ and the graphical
variables $K_{n}$ (of which finitely many will appear for 
any given virtual knot or link). ${\cal A}[K]$ is a regular isotopy invariant of virtual knots and links. The normalized version
$${\cal W}[K] = (-A^{3})^{-wr(K)} {\cal A}[K]$$ is an invariant virtual isotopy. If we set $A = 1$ and $d = -A^2 - A^{-2} = -2$, then the resulting specialization
$${\cal F}[K] = {\cal A}[K](A = 1)$$
is an invariant of flat virtual knots and links.
\bigbreak

\noindent {\bf Proof Sketch.} We have already discussed the first Reidemeister move. Invariance under virtual detours is implicit in the definition of the 
state sum. The arrow polynomial has framing behaviour under the first Reidemeister moves as shown in Figure 18.
In Figures 19 and 20 we show how the removal of a disoriented loop gives rise to invariance under the directly oriented
second Reidemeister move. In the state expansion, the disoriented loop is removed but multiplies its term by $d = -A^{2} - A^{-2}.$  The other two disoriented local configurations receive vertex weights
of $A^2$ and $A^{-2}$, and {\it each of these configurations has the same set of reduced states}. 
Thus these three configurations cancel each other from the state sum. This leaves the remaining local state with parallel arcs, 
and gives invariance under the directly oriented second Reidemeister move. Invariance under the  reverse oriented second Reidemeister move is shown in Figure 20. Invariance under the third Reidemeister move follows from an analysis of Figure 21.
\bigbreak

\begin{figure}
     \begin{center}
     \begin{tabular}{c}
     \includegraphics[width=7cm]{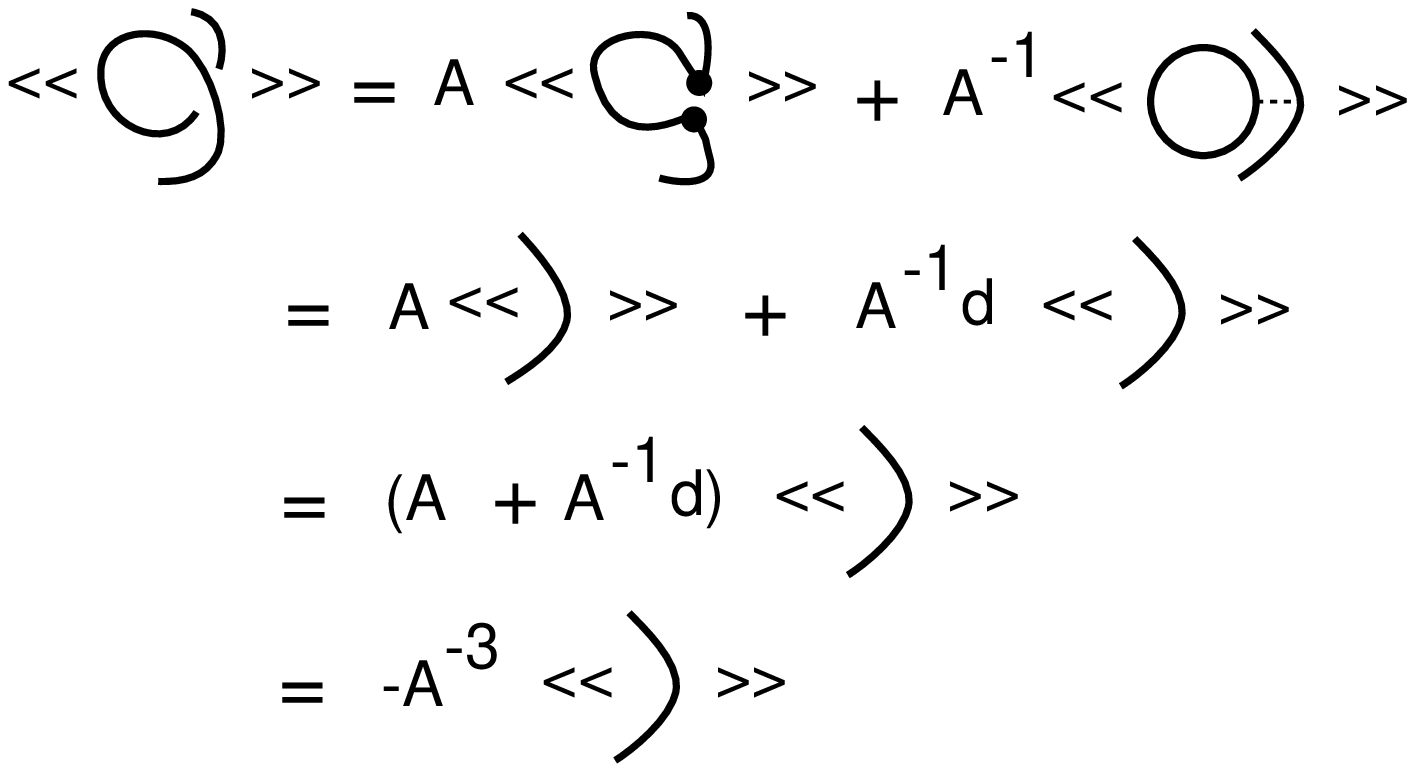}
     \end{tabular}
     \caption{\bf The Type One Move}
     \label{Figure 18}
\end{center}
\end{figure}

\begin{figure}
     \begin{center}
     \begin{tabular}{c}
     \includegraphics[width=9cm]{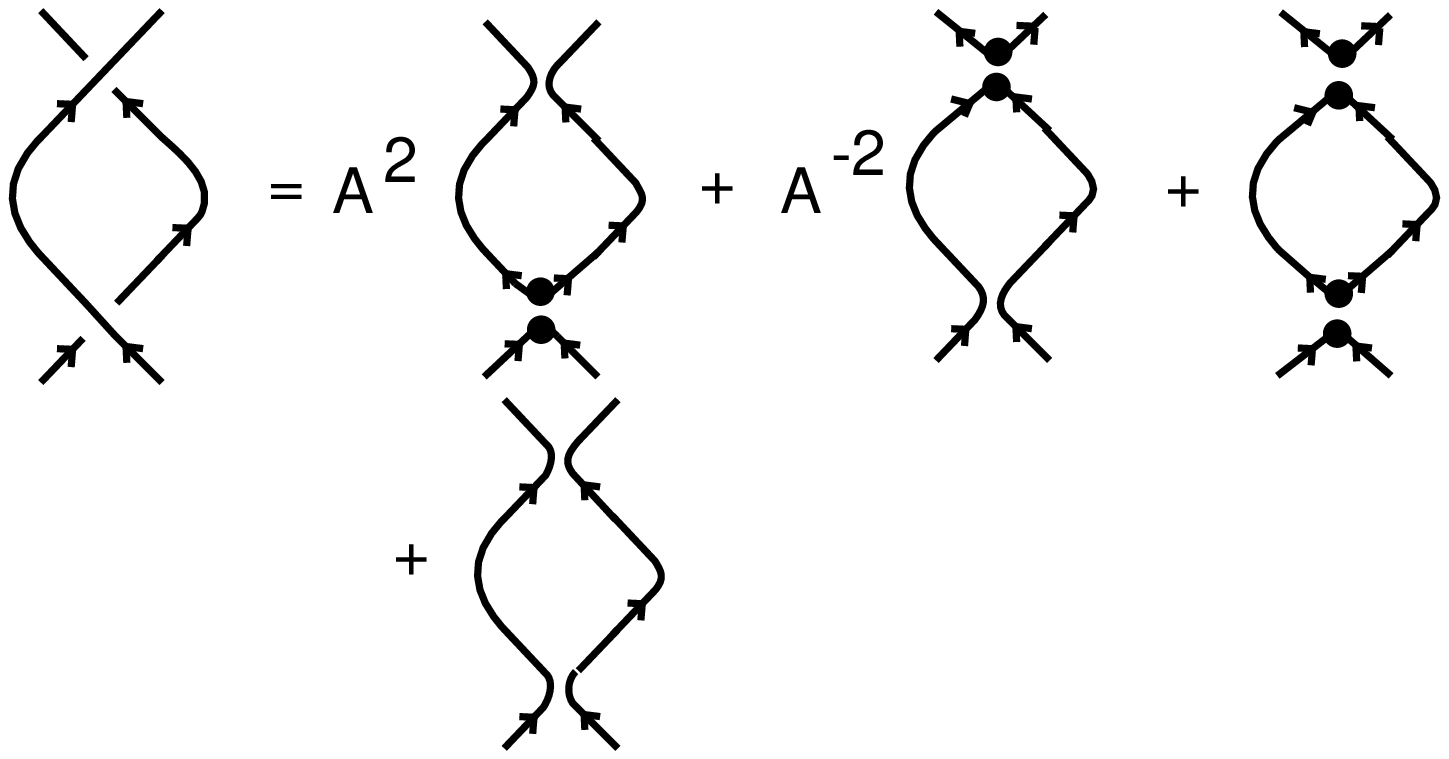}
     \end{tabular}
     \caption{\bf Oriented Second Reidemeister Move}
     \label{Figure 19}
\end{center}
\end{figure}

\begin{figure}
     \begin{center}
     \begin{tabular}{c}
     \includegraphics[width=8cm]{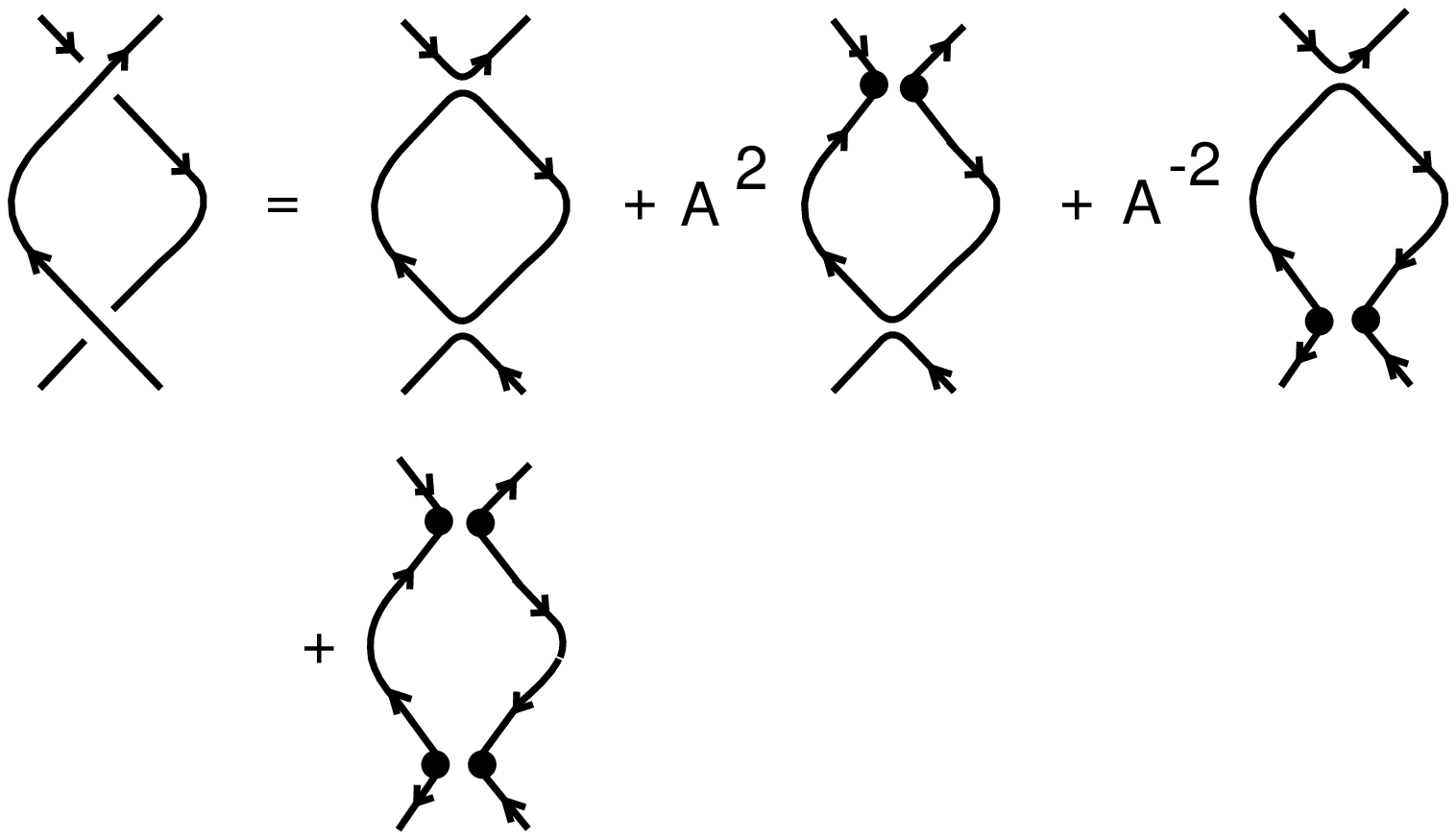}
     \end{tabular}
     \caption{\bf Reverse Oriented Second Reidemeister Move}
     \label{Figure 20}
\end{center}
\end{figure}

\begin{figure}
     \begin{center}
     \begin{tabular}{c}
     \includegraphics[width=10cm]{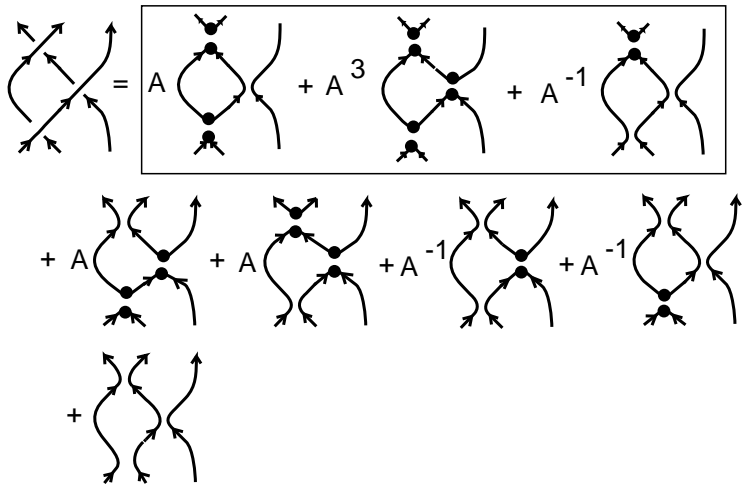}
     \end{tabular}
     \caption{\bf Third Reidemeister Move}
     \label{Figure 21}
\end{center}
\end{figure}

Here is a first example of a calculation of the arrow polynomial invariant. View Figure 22. The virtual knot $K$ in this figure has two crossings.
One can see that this knot is a non-trival virtual knot by simply calculating the odd writhe $J(K)$ (defined in section 5). We have that
$J(K) = 2,$ proving that $K$ is non-trivial and non-classical. This is the simplest virtual knot, the analog of the trefoil knot for virtual knot 
theory. The arrow polynomial polynomial gives an independent verification that $K$ is non-trivial and non-classical.
\bigbreak

\begin{figure}
     \begin{center}
     \begin{tabular}{c}
     \includegraphics[width=10cm]{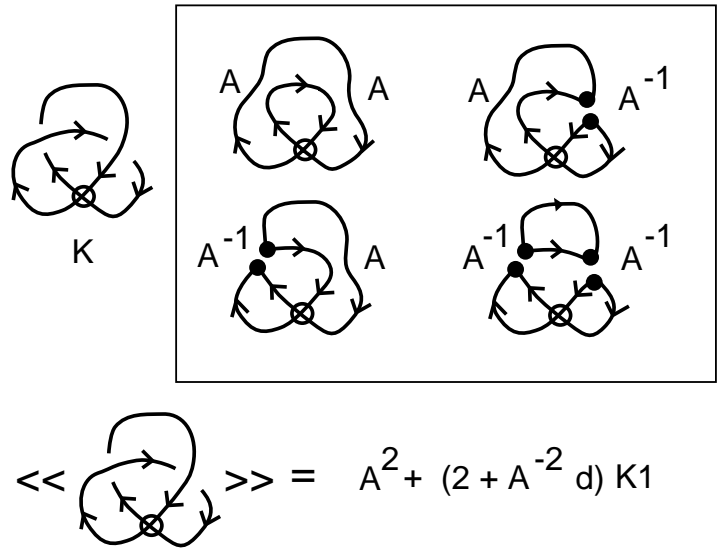}
     \end{tabular}
     \caption{\bf Example1}
     \label{Figure 22}
\end{center}
\end{figure}

In the next example, shown in Figure 23, we have a long virtual diagram $L$ with two crossings. The calculation of the arrow polynomial for 
$L$ is given in Figure 23 and shows that it is a non-trivial and non-classical. In fact, this same formalism
proves that $Flat(L)$ is a non-trivial flat link (as we have also by using the odd writhe.). {\it Note that the arrow polynomial is an invariant of 
flat diagrams when we take $A=1.$} In the Figure 23 we have illustrated for this example how 
non-cancelling cusps can occur on a long state in the expasion of the invariant. This means that, 
effectively, there are more variables for the arrow polynomial as an invariant of long knots, and this
gives it extra powers of discrimination for in the long category. Note the closure of $L$ of Figure 23 is a trivial knot. Thus we see that the arrow polynomial can sometimes discriminate between a long knot and its closure. Another example of this type is shown
in Figure 24. In both of these cases the cusps on the long state are in an order that would be differerent
if we were to turn the long knot by $180$ degrees. Thus the invariant distinguishes between the 
long knot and its reverse.
\bigbreak

\begin{figure}
     \begin{center}
     \begin{tabular}{c}
     \includegraphics[width=10cm]{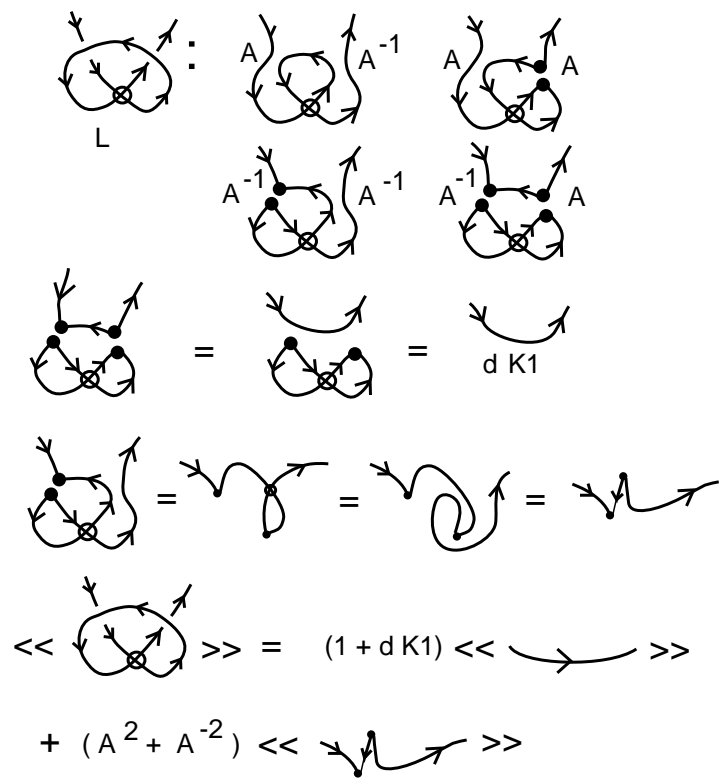}
     \end{tabular}
     \caption{\bf Example2}
     \label{Figure 23}
\end{center}
\end{figure}

\begin{figure}
     \begin{center}
     \begin{tabular}{c}
     \includegraphics[width=10cm]{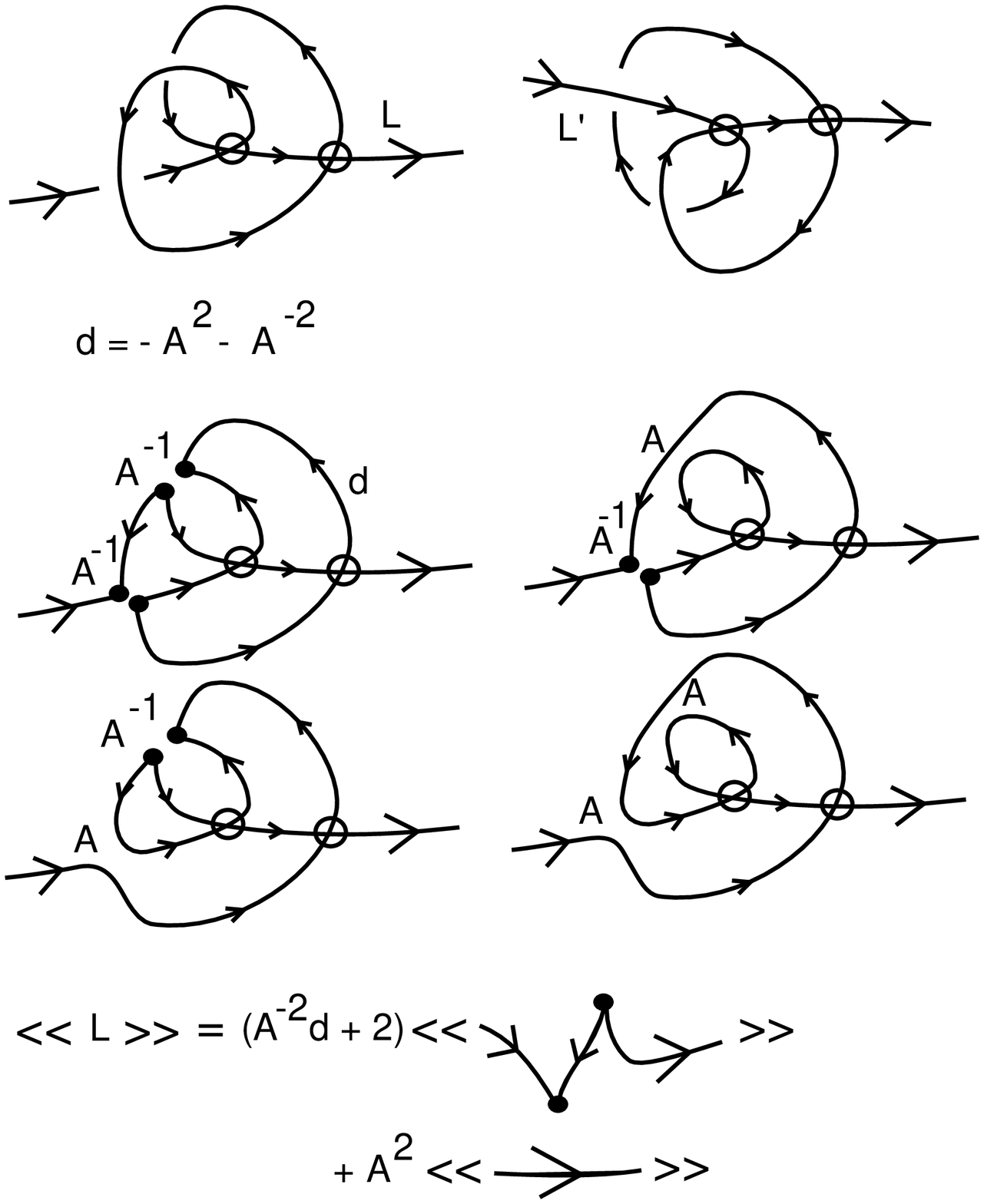}
     \end{tabular}
     \caption{\bf Example3}
     \label{Figure 24}
\end{center}
\end{figure}

The next example is given in Figure 25. Here we calculate the arrow polynomial for a non-trivial virtual knot with unit Jones polynomial.  Specialization of the calculation to $A=1$ shows that the corresponding flat knot is non-trivial as well. 
\bigbreak

Figures 26 exhibits the calculation of the arrow polynomial for the Kishino diagram, showing once 
again that the Kishino knot is non-trivial and that its underlying flat diagram is also non-trivial.
\bigbreak

The example shown in Figures 27 and Figure 28 shows the result of expanding a virtualized
classical crossing using the arrow polynomial state sum. Virtualization of a crossing was described in Section 6. In a
virtualized crossing, one sees a classical crossing that is flanked by two virtual crossings. In Section 6  we showed that
the standard bracket state sum does not see virtualization in the sense that it  has the same value as the result of
smoothing both flanking virtual crossings that have been added to the diagram. The result is that the value of of the
bracket polynomial of the knot with a virtualized classical crossing is the same as the value of the bracket polynomial of
the original knot after the same crossing has been {\it switched} (exchanging over and undercrossing segments). 
\bigbreak

As one can see from the formula in Figure 27, this smoothing property of the bracket polynomial will not generally be the case for the arrow polynomial state
sum. In Figure 28 we show that this difference is indeed the case for an infinite collection of examples. In that figure we use a tangle $T$ that is assumed
to be a classical tangle. arrow polynomial expansison of this tangle is necessarily of the form shown in that figure: a linear combination of an oriented
smoothing and a reverse oriented smoothing with respective coefficients $a(T)$ and $b(T)$ in the Laurent polynomial ring $Q[A,A^{-1}].$ We leave the
verification of this fact to the reader. In Figure 28 we show a generic diagram that is obtained by a {\it single} virtualization from a classical diagram,
and we illustrate the calculation of its arrow polynomial invariant. As the reader can see from this Figure, there is a non-trivial graphical term whenever
$b(T)$ is non-zero. Thus we conclude that the single virtualization of any classical link diagram (in the form shown in this figure) will be non-trivial and
non-classical whenever $b(T)$ in non-zero. This is an infinite class of examples, and the result can be used to recover the results about single
virtualization that we  obtained in a previous paper with Heather Dye \cite{MinSurf} using the surface bracket polynomial.
\bigbreak

For more information about the arrow polynomial, we refer the reader to our paper \cite{DyeKauff} where
we prove that the maximal monomial degree in the $K_{n}$ variables  is a lower bound for the virtual crossing number of the virtual knot or link. There are many open problems associated with this 
estimate for the virtual crossing number. Also the reader of \cite{DyeKauff} will encounter examples of virtual knots and links that are undetectable by the arrow polynomial.
\bigbreak

\begin{figure}
     \begin{center}
     \begin{tabular}{c}
     \includegraphics[width=7cm]{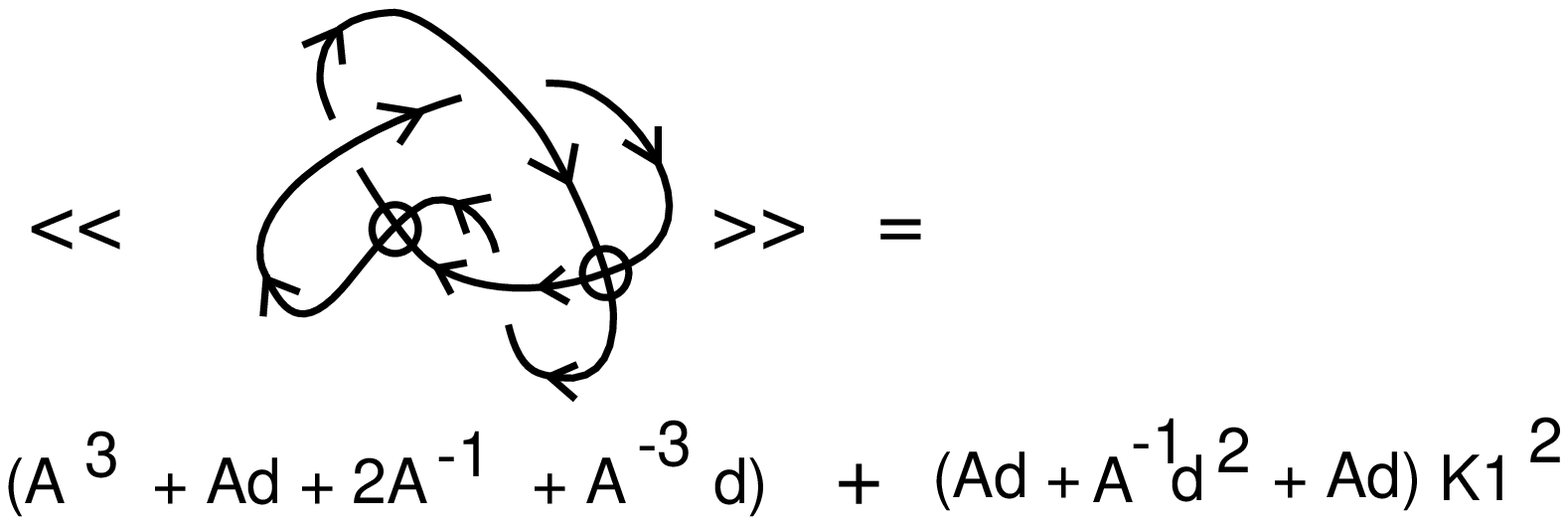}
     \end{tabular}
     \caption{\bf Arrow Polynomial for the Virtualized Trefoil}
     \label{Figure 25}
\end{center}
\end{figure}

\begin{figure}
     \begin{center}
     \begin{tabular}{c}
     \includegraphics[width=8cm]{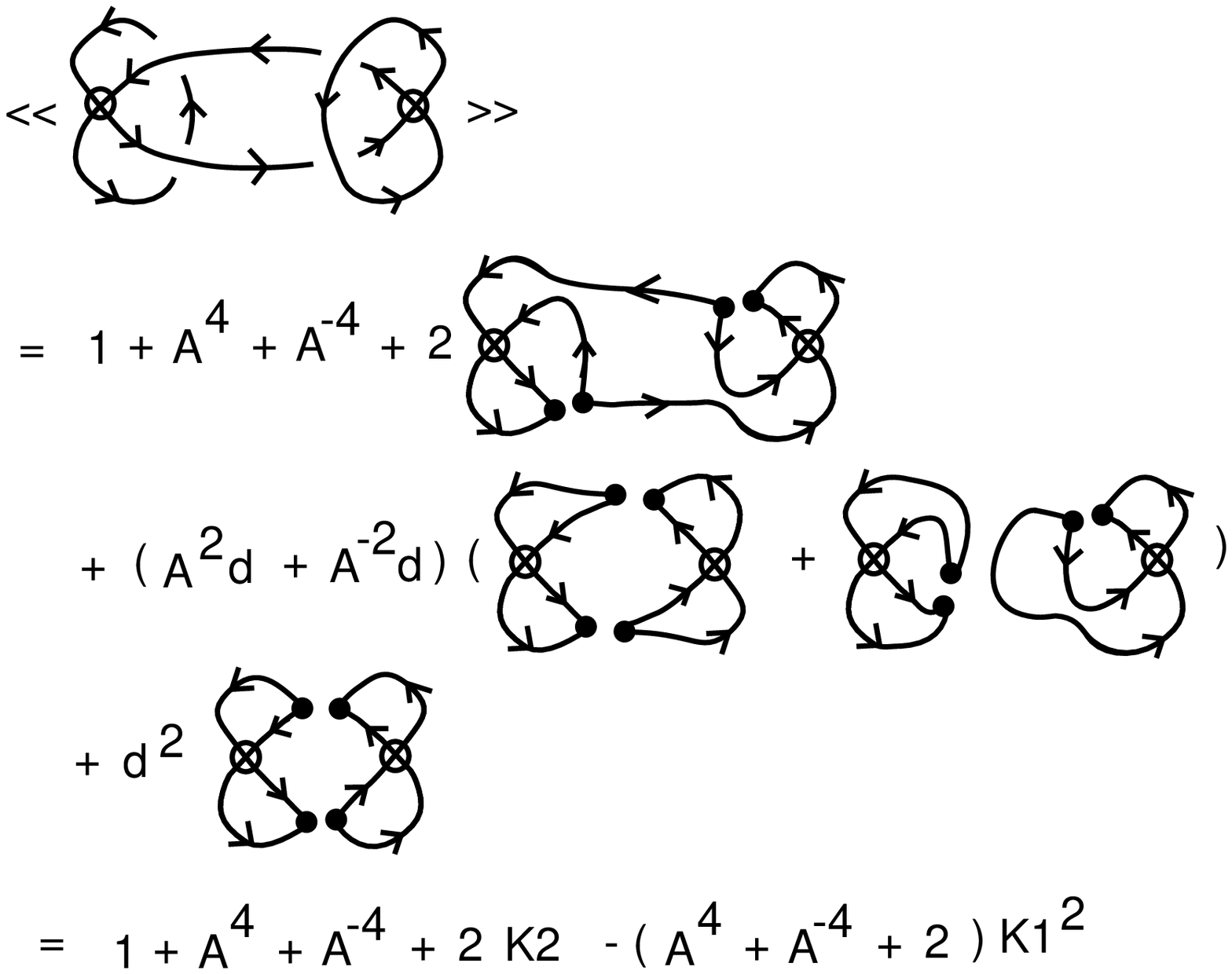}
     \end{tabular}
     \caption{\bf Arrow Polynomial for the Kishino Diagram}
     \label{Figure 26}
\end{center}
\end{figure}

\begin{figure}
     \begin{center}
     \begin{tabular}{c}
     \includegraphics[width=8cm]{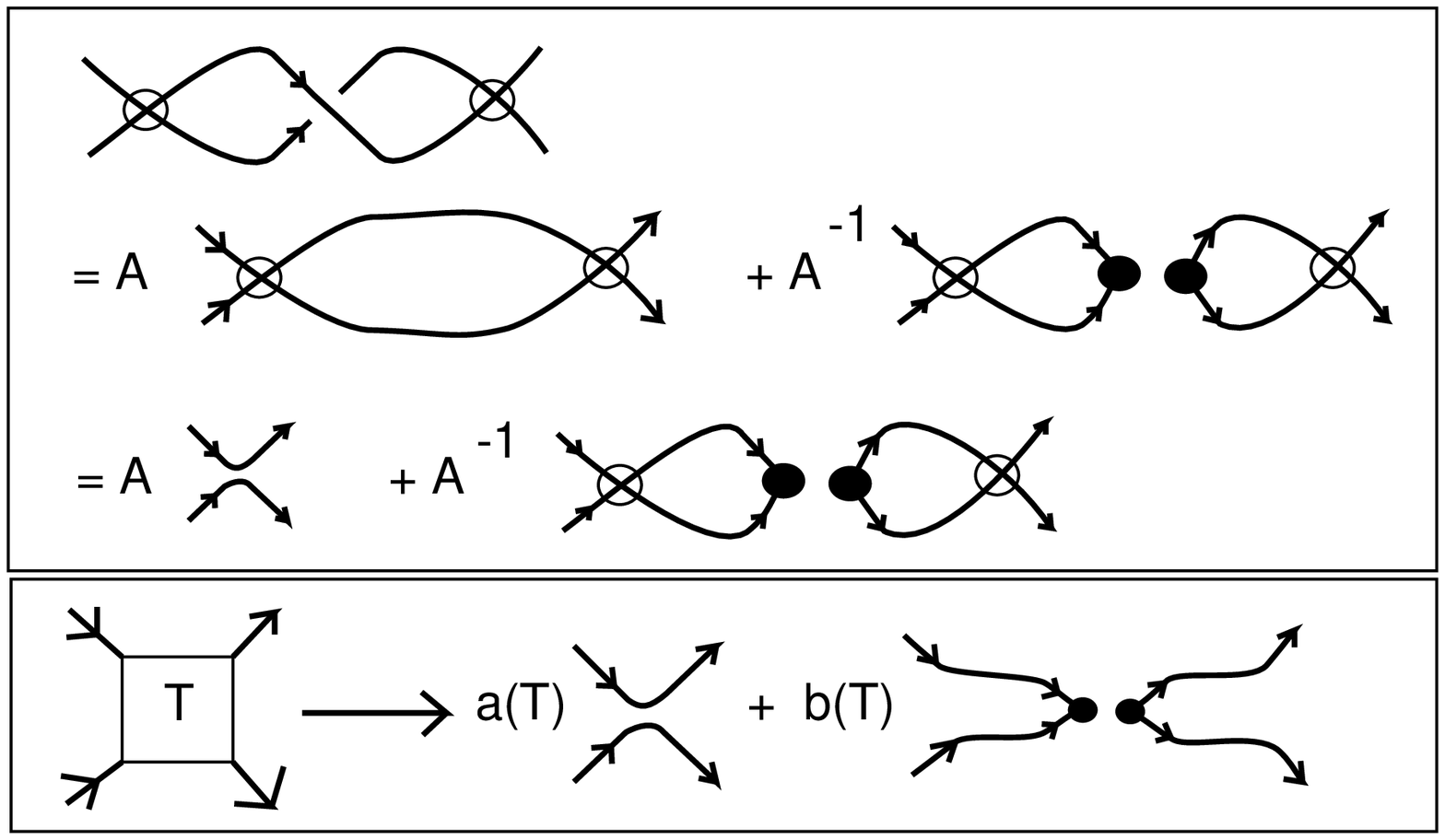}
     \end{tabular}
     \caption{\bf The $1$-Virtualization of a Classical Diagram (A)}
     \label{Figure 27}
\end{center}
\end{figure}

\begin{figure}
     \begin{center}
     \begin{tabular}{c}
     \includegraphics[width=10cm]{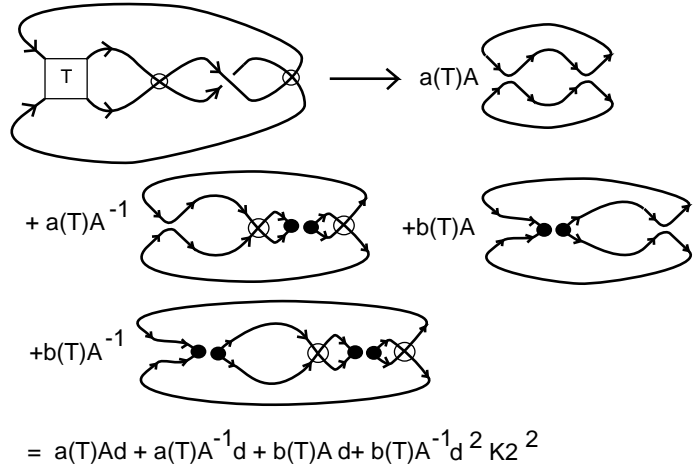}
     \end{tabular}
     \caption{\bf The $1$-Virtualization of a Classical Diagram (B)}
     \label{Figure 28}
\end{center}
\end{figure}
\bigbreak
 
 \subsection{Generalizing Parity to the Arrow Polynomial}
The purpose of this section is to indicate how our construction of the Parity Bracket in Section 7
generalizes to the arrow polynomial. The reader will find it easy to verify that at the site of a Reidemeister
Two move either both crossings are odd, or both are even. At the site of a Reidemeister Three move
either all three crossings are even, or two of the three are odd. It follows from this that if we decide
(as in the Parity Bracket) to make graphical vertices from all the odd crossings and expand the 
oriented bracket (for the arrow polynomial) then we must factor the graphs that result from the state expansion on the even crossings of the diagram by an equivalence relation generated by patterns of 
graph reduction corresponding to the second Reidemeister moves on graphical vertices. These reductions come in a number of possible orientations, and it is possible to have a single cusp
involved as illustrated in Figure 29. In this figure we do not indicate the orientations. Because the
cusps in the graphs are involved with the graphical reductions one takes the resulting graphs up to 
virtual graphical equivalence, these reductions and the cusp cancellation is restricted to single edges of the graphs. The result is a parity version of the arrow polynomial. We will not calculate examples of this invariant in this survey. The Parity Arrow Polynomial will be taken up in a subsequent paper.
\bigbreak

\begin{figure}
     \begin{center}
     \begin{tabular}{c}
     \includegraphics[width=10cm]{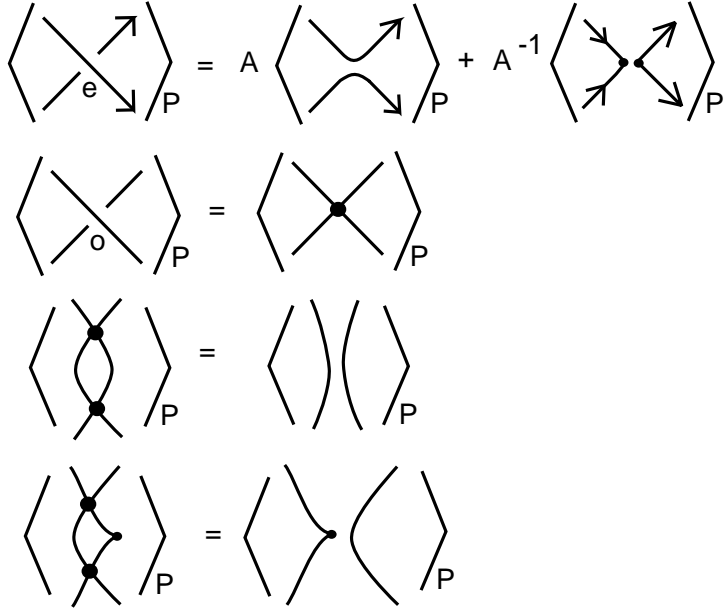}
     \end{tabular}
     \caption{\bf Relations for the Parity Arrow Polynomial}
     \label{Figure 29}
\end{center}
\end{figure}
\bigbreak

\section{Khovanov Homology}

In this section, we describe Khovanov homology
along the lines of \cite{Kho,DB1}, and we tell the story so that the gradings and the structure of the differential emerge in a natural way.
This approach to motivating the Khovanov homology uses elements of Khovanov's original approach, Viro's use of enhanced states for the bracket
polynomial \cite{V}, and Bar-Natan's emphasis on tangle cobordisms \cite{DB2}. We use similar considerations in our paper \cite{DKM}. The purpose of this section is to provide background for the 
next section where we describe a categorifcation of the arrow polynomial and give examples of
detection that can be accomplished by this homological invariant.
\bigbreak

Two key motivating ideas are involved in finding the Khovanov invariant. First of all, one would like to {\it categorify} a link polynomial such as
$\langle K \rangle.$ There are many meanings to the term categorify, but here the quest is to find a way to express the link polynomial
as a {\it graded Euler characteristic} $\langle K \rangle = \chi_{q} \langle H(K) \rangle$ for some homology theory associated with $\langle K \rangle.$
\bigbreak

The bracket polynomial \cite{K,KNOTS} model for the Jones polynomial \cite{Jones,JO1,JO2,WITT} is usually described by the expansion
\begin{equation}
\langle \Across \rangle=A \langle \Asmooth \rangle + A^{-1}\langle
\Bsmooth \rangle \label{kabr}
\end{equation}

and we have

\begin{equation}
\langle K \, \bigcirc \rangle=(-A^{2} -A^{-2}) \langle K \rangle \label{kabr}
\end{equation}

\begin{equation}
\langle \Rcurl \rangle=(-A^{3}) \langle \Arc \rangle \label{kabr}
\end{equation}

\begin{equation}
\langle \Lcurl \rangle=(-A^{-3}) \langle \Arc \rangle \label{kabr}
\end{equation}
\bigbreak

Letting $c(K)$ denote the number of crossings in the diagram $K,$ if we replace $\langle K \rangle$ by 
$A^{-c(K)} \langle K \rangle,$ and then replace $A$ by $-q^{-1},$ the bracket will be rewritten in the
following form:
\begin{equation}
\langle \Across \rangle=\langle \Asmooth \rangle-q\langle
\Bsmooth \rangle \label{kabr}
\end{equation}
with $\langle \bigcirc\rangle=(q+q^{-1})$.
It is useful to use this form of the bracket state sum
for the sake of the grading in the Khovanov homology (to be described below). We shall
continue to refer to the smoothings labeled $q$ (or $A^{-1}$ in the
original bracket formulation) as {\it $B$-smoothings}. We should
further note that we use the well-known convention of {\it enhanced
states} where an enhanced state has a label of $1$ or $X$ on each of
its component loops. We then regard the value of the loop $q + q^{-1}$ as
the sum of the value of a circle labeled with a $1$ (the value is
$q$) added to the value of a circle labeled with an $X$ (the value
is $q^{-1}).$ We could have chosen the more neutral labels of $+1$ and $-1$ so that
$$q^{+1} \Longleftrightarrow +1 \Longleftrightarrow 1$$
and
$$q^{-1} \Longleftrightarrow -1 \Longleftrightarrow X,$$
but, since an algebra involving $1$ and $X$ naturally appears later, we take this form of labeling from the beginning.
\bigbreak

To see how the Khovanov grading arises, consider the form of the expansion of this version of the 
bracket polynonmial in enhanced states. We have the formula as a sum over enhanced states $s:$
$$\langle K \rangle = \sum_{s} (-1)^{n_{B}(s)} q^{j(s)}$$
where $n_{B}(s)$ is the number of $B$-type smoothings in $s$, $\lambda(s)$ is the number of loops in $s$ labeled $1$ minus the number of loops
labeled $X,$ and $j(s) = n_{B}(s) + \lambda(s)$.
This can be rewritten in the following form:
$$\langle K \rangle  =  \sum_{i \,,j} (-1)^{i} q^{j} dim({\cal C}^{ij})$$
where we define ${\cal C}^{ij}$ to be the linear span (over $k = Z/2Z$ as we will work with mod $2$ coefficients) of the set of enhanced states with $n_{B}(s) = i$ and $j(s) = j.$
Then the number of such states is the dimension $dim({\cal C}^{ij}).$ 
\bigbreak

\noindent We would like to have a  bigraded complex composed of the ${\cal C}^{ij}$ with a
differential
$$\partial:{\cal C}^{ij} \longrightarrow {\cal C}^{i+1 \, j}.$$ 
The differential should increase the {\it homological grading} $i$ by $1$ and preserve the 
{\it quantum grading} $j.$
Then we could write
$$\langle K \rangle = \sum_{j} q^{j} \sum_{i} (-1)^{i} dim({\cal C}^{ij}) = \sum_{j} q^{j} \chi({\cal C}^{\bullet \, j}),$$
where $\chi({\cal C}^{\bullet \, j})$ is the Euler characteristic of the subcomplex ${\cal C}^{\bullet \, j}$ for a fixed value of $j.$
\bigbreak

\noindent This formula would constitute a categorification of the bracket polynomial. Below, we
shall see how {\it the original Khovanov differential $\partial$ is uniquely determined by the restriction that $j(\partial s) = j(s)$ for each enhanced state
$s$.} Since $j$ is 
preserved by the differential, these subcomplexes ${\cal C}^{\bullet \, j}$ have their own Euler characteristics and homology. We have
$$\chi(H({\cal C}^{\bullet \, j})) = \chi({\cal C}^{\bullet \, j}) $$ where $H({\cal C}^{\bullet \, j})$ denotes the homology of the complex 
${\cal C}^{\bullet \, j}$. We can write
$$\langle K \rangle = \sum_{j} q^{j} \chi(H({\cal C}^{\bullet \, j})).$$
The last formula expresses the bracket polynomial as a {\it graded Euler characteristic} of a homology theory associated with the enhanced states
of the bracket state summation. This is the categorification of the bracket polynomial. Khovanov proves that this homology theory is an invariant
of knots and links (via the Reidemeister moves) , creating a new and stronger invariant than the original Jones polynomial.
\bigbreak

The differential is based on regarding two states as {\it adjacent} if one differs from the other by a single smoothing at some site.
Thus if $(s,\tau)$ denotes a pair consisting in an enhanced state $s$ and site $\tau$ of that state with $\tau$ of type $A$, then we consider
all enhanced states $s'$ obtained from $s$ by smoothing at $\tau$ and relabeling only those loops that are affected by the resmoothing.
Call this set of enhanced states $S'[s,\tau].$ Then we shall define the {\it partial differential} $\partial_{\tau}(s)$ as a sum over certain elements in
$S'[s,\tau],$ and the differential by the formula $$\partial(s) = \sum_{\tau} \partial_{\tau}(s)$$ with the sum over all type $A$ sites $\tau$ in $s.$
It then remains to see what are the possibilities for $\partial_{\tau}(s)$ so that $j(s)$ is preserved.
\bigbreak

\noindent {\bf Proposition.} The partial differentials $\partial_{\tau}(s)$ are uniquely determined by the condition that $j(s') = j(s)$ for all $s'$
involved in the action of the partial differential on the enhanced state $s.$ This unique form of the partial differential can be described by the 
following structures of multiplication and comultiplication on the algebra \cal{A} = $k[X]/(X^{2})$ where $k = Z/2Z$ for mod-2 coefficients, or $k = Z$
for integral coefficients.
\begin{enumerate}
\item The element $1$ is a multiplicative unit and $X^2 = 0.$
\item $\Delta(1) = 1 \otimes X + X \otimes 1$ and $\Delta(X) = X \otimes X.$
\end{enumerate}
These rules describe the local relabeling process for loops in a state. Multiplication corresponds to the case where two loops merge to a single loop, 
while comultiplication corresponds to the case where one loop bifurcates into two loops.
(The proof is omitted.)
\bigbreak 

There is more to say about the nature of this construction with respect to Frobenius algebras and tangle
cobordisms. The partial boundaries can be conceptualized in terms of surface cobordisms. The equality of mixed partials corresponds to topological equivalence of the corresponding surface cobordisms, and to the relationships between Frobenius algebras and the surface
cobordism category. The proof of invariance of Khovanov homology with respect to the Reidemeister moves (respecting grading changes) will not be given here.
See \cite{Kho,DB1,DB2}. It is remarkable that this version of Khovanov homology is uniquely specified by natural ideas about adjacency of states in the bracket
polynomial.
\bigbreak

\section{Categorifying  the Arrow Polynomial $\langle \langle K \rangle \rangle = {\cal A}[K]$}
It is not hard to see that Khovanov homology modulo two extends to a homology theory for 
virtual links. In \cite{DKM} we extend this to a new homology theory  that categorifies the arrow polynomial.
\bigbreak

In order to consider gradings for Khovanov homology in relation to
the arrow polynomial ${\cal A}[K]$ we have to examine how the
arrow number of state loops change under a replacement of an
$A$-smoothing by a $B$-smoothing. Such replacement, when we use
oriented diagrams involves the replacement of a cusp pair by an
oriented smoothing or vice versa. Furthermore, we may be combining
or splitting two loops. Refer to Figure 30  for a
depiction of the different cases. This figure shows the three basic
cases. 
These are all the ways that loops can interact and change their respective arrow numbers. We apply these results to the grading in Khovanov homology.
\bigbreak

\begin{figure}
     \begin{center}
     \begin{tabular}{c}
     \includegraphics[width=7cm]{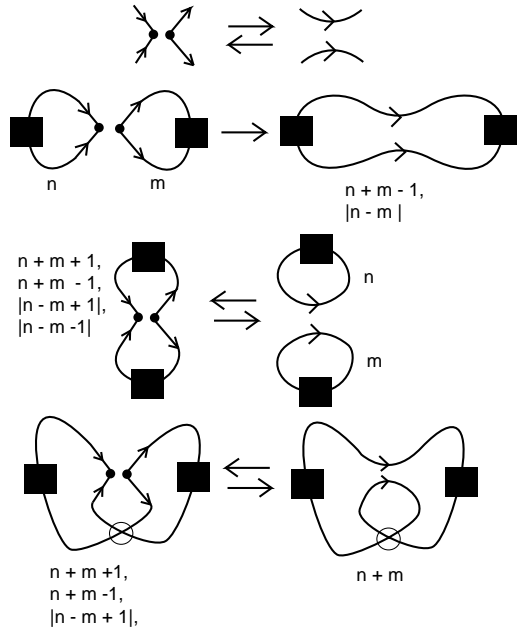}
     \end{tabular}
     \caption{\bf Arrow Numbers for Interacting Loops}
     \label{Figure 30}
\end{center}
\end{figure}

\subsection{Dotted gradings and the dotted categorification}
We call an enhanced state {\it dotted} if it has an odd arrow number, and otherwise it is undotted.
The number of dotted $X$'s and the number of dotted $1$'s in the state correspond to the number
of enhanced state loops (labeled with $X$ or with $1$) that are dotted. We define a grading
$$g(s) = \sharp ({\dot X}) - \sharp ({\dot 1})$$
\bigbreak

\noindent {\bf Theorem.}
 Take $[[K]]_{g}$ to be the space of enhanced states for $K$ endowed with new grading
as above.
Define $\partial'$ to be the composition of $\partial$ with the new grading projection
and set $\partial''=\partial-\partial'$.
Then the homology of $[[K]]_{g}$ (with respect to $\partial'$) is invariant (up to a degree shift
and a height shift).
\bigbreak

With this result we have a version of Khovanov homology that uses the arrow information.
In \cite{DKM} we have investigated this categorification and other related ones for the arrow polynomial.
So far this  categorifcation is effective in that {\it we \cite{KaestKauff} have many pairs of 
virtual knots that are distinguished by the arrow categorification, but not distinguished by either
Khovanov homology or the arrow polynomial.} An example of such a pair (due to Aaron Kaestner) is shown in Figure 31 and 32. 
In Figure 31 we illustrate the knot $VK5[129]$ and in Figure 32 we illustrate the knot $VK5[267]$ from 
Jeremy Green's tables (See the website of Dror Bar Natan for these tables.). The examples shown in these figures do not seem to be related by any form of mutation and the difference betweem them is not dectected by the arrow polynomial, but it is detected by the categorification of the arrow polynomial.
\bigbreak

Note however, that one can compute the Parity Bracket (See Section 7) for these examples. In the case of $VK5[129]$ there are five crossings and four odd crossings. It is easy to see that the Parity Bracket expansion contains an irreducible graph with four nodes. On the other hand, the diagram for
$VK5[267]$  has two odd crossings and hence its Parity Bracket (which one can compute to be non-trivial with an irreducible graphs with two nodes) cannot be equal up to normalization with the Parity
Bracket for $VK5[129].$ Thus this pair of knots is distinguished by using parity. This example shows
the power of the parity method.
\bigbreak

This approach through examples and combining categorification and parity is very fruitful for producing examples and for sharpening concepts. We are programming an integral version of Khovanov homology for virtual knots that is structured along the lines proposed by Vassily Manturov \cite{M}, and we expect all this work to reflect back on better understanding of Khovanov homology for classical knots. This is the first instance of non-trivial use of categorified homology in virtual knot theory.

\begin{figure}
     \begin{center}
     \begin{tabular}{c}
     \includegraphics[width=6cm]{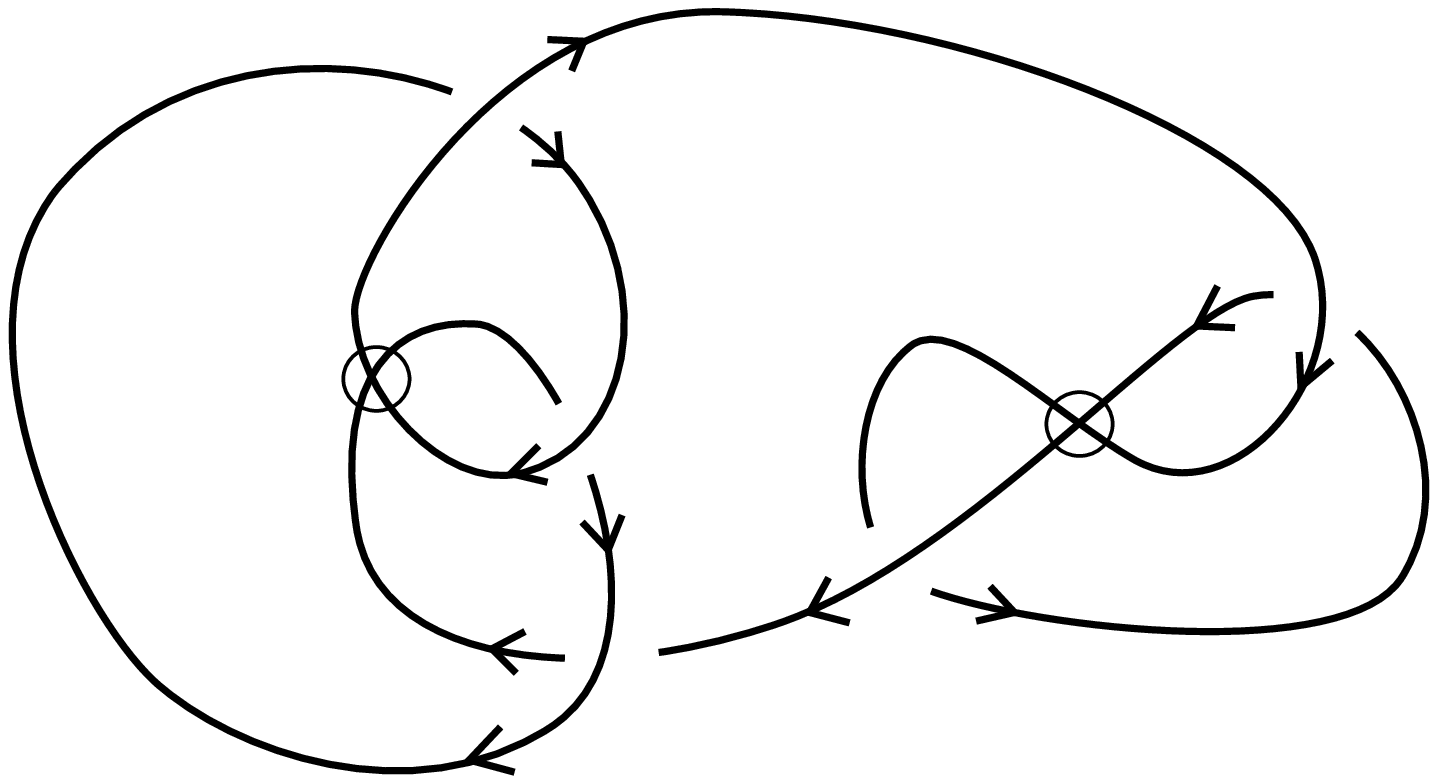}
     \end{tabular}
     \caption{\bf VK5[129]}
     \label{Figure 31}
\end{center}
\end{figure}

\begin{figure}
     \begin{center}
     \begin{tabular}{c}
     \includegraphics[width=6cm]{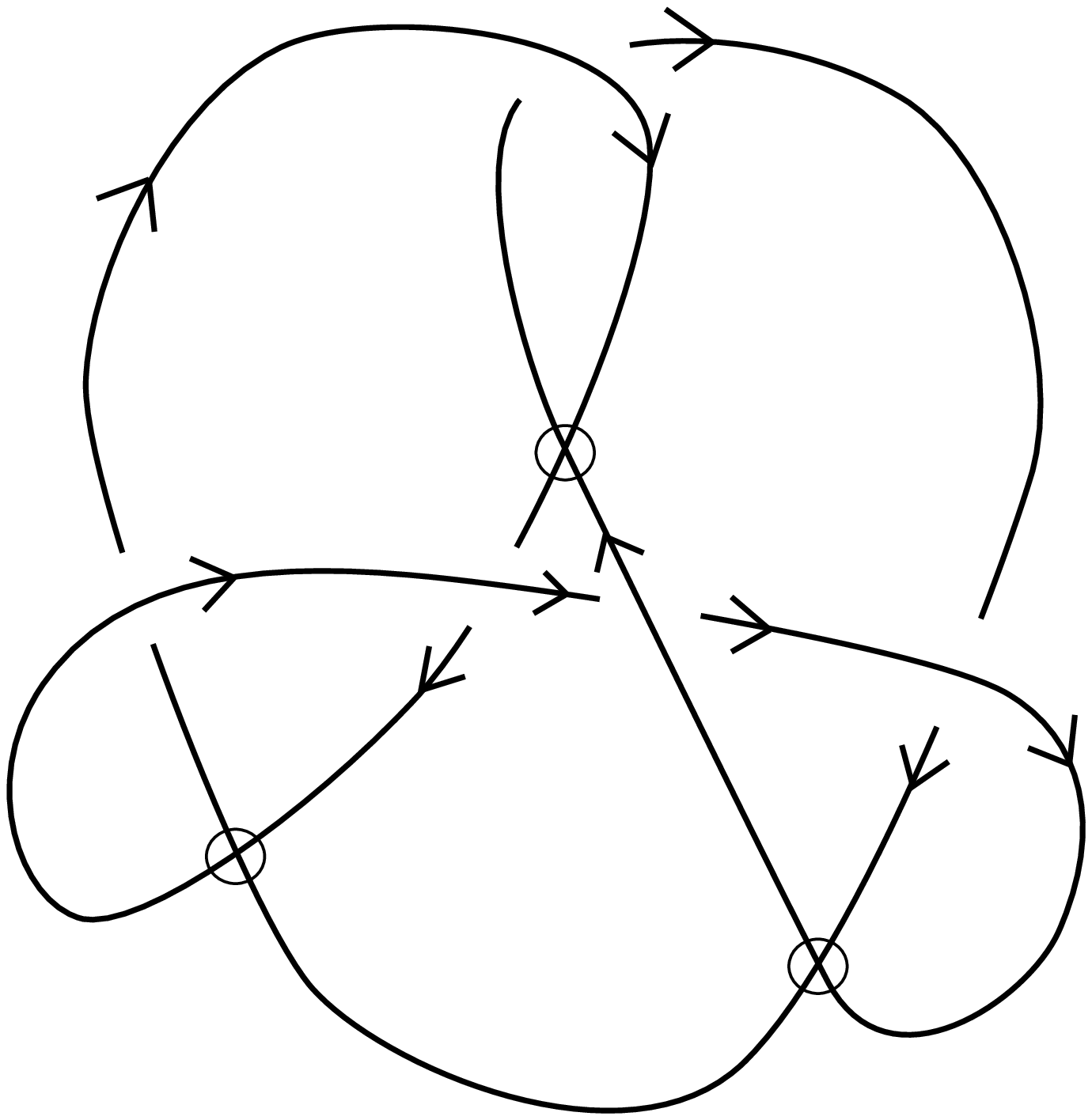}
     \end{tabular}
     \caption{\bf VK5[267]}
     \label{Figure 32}
\end{center}
\end{figure}

\section {Discussion} 
In this paper we have given an introduction to the subject of virtual knot theory and we have concentrated on the easy extension of the Jones polynomial to virtual knots and links and a more subtle 
generalization that we have called the {\it arrow polynomial}. The arrow polynomial uses the oriented structure of link diagrams in a way that is particularly suited to virtual knot theory and obtains an 
invariiant of virtual knots and links that has infinitely many variables. The ideas behind the arrow polynomial have been used to create an extended bracket polynomial \cite{ExtBr} that uses
more of the oriented structure and has graphs representing its extra variables. The arrow polynomial itself can be seen to be a reformulation or generalization of some of the constructions of Miyazawa and Kamada \cite{Miyazawa1,Miyazawa2,Kamada1,Kamada2,KamadaMiya}. After discussing the 
arrow polynomial, we gave a quick summary of Khovanov homology and a brief description of our categorification  of the arrow polynomial \cite{DyeKauff,DKM} and an indication of the fact that this categorification gives
a yet stronger invariant and generates even more new problems. 
\bigbreak

We also designed this paper to present examples of the use of parity in virtual knot theory. We began with a discussion of the odd writhe in Section 5.   In the Gauss code of the diagram each crossing appears twice with a sequence of symbols in between the two appearances. When this sequence is odd we say that the given crossing is odd. When this sequence is even, we say that the given crossing is even. Many possibilites for new invariants are inherent in this parity distinction between even and odd crossings. We have given one example of this kind of application in Section 7 where we discuss the Parity Bracket  of Vassily Manturov and apply it to solving a conjecture that we had made about virtual knot diagrams with unit  Jones polynomial. A second example is given at the end of the paper to distinguish the knots in Figures 31 and 32. There is much potential for the use of parity in virtual knot theory.
\bigbreak

\end{document}